 \documentclass[12pt]{article}
  \usepackage{amsfonts,amsmath,amssymb}
  \usepackage[square,sort&compress,comma,numbers]{natbib}
   \def\R{\mathbb{R}}
   \def\N{\mathbb{N}}
   \def\Z{\mathbb{Z}}
   \def\1{{\rm I\mskip -10.5mu 1}}
   
   \def\d{{\delta}}
   \def\e{{\varepsilon}}
   \def\D{{\nabla}}

   \def\vi{{\varphi}}

   \def\cA{{\cal A}}
   
   \def\cC{{\cal C}}
   \def\cD{{\cal D}}
   
   \def\cG{{\cal G}}
   \def\cH{{\cal H}}

   \def\cN{{\cal N}}

   \def\dist{\mathop{\rm dist}\nolimits}
   
   \def\const{\mathop{\rm const}\nolimits}
   
   \def\loc{\mathop{\rm loc}\nolimits}

   \def\rad{\mathop{\rm rad}\nolimits}
   \def\no{\noindent}
   \def\proof{\noindent{\bf Proof \hspace{2mm}}}
   \def\qed{{\hfill {\em q.e.d.}\\\vspace{1mm}}}

   \newcommand{\beq}{\begin{equation}}
   \newcommand{\eeq}{\end{equation}}

\newtheorem{df}{Definition}[section]
\newtheorem{prop}[df]{Proposition}
\newtheorem{lemma}[df]{Lemma}
\newtheorem{teo}[df]{Theorem}
\newtheorem{rem}[df]{Remark}

\newtheorem{cor}[df]{Corollary}

 \newcommand{\sezione}[1]{\section{#1}\setcounter{equation}{0}}

\begin{document}


  \title{Positive solutions for autonomous and non-autonomous
     nonlinear critical elliptic problems in exterior domains
}
 \vspace{5mm}


  \maketitle


\begin{center}

{ {\bf Sergio LANCELOTTI$^a$,\quad  Riccardo MOLLE$^b$}}

\vspace{5mm}

{\em
${\phantom{1}}^a$Dipartimento di Scienze Matematiche,
Politecnico di Torino,\linebreak
Corso Duca degli Abruzzi n. 24, 10129 Torino, Italy.\linebreak
sergio.lancelotti@polito.it
}

\vspace{2mm}

{\em
${\phantom{1}}^b$Dipartimento di Matematica,
Universit\`a di Roma ``Tor Vergata'',\linebreak
Via della Ricerca Scientifica n. 1,
00133 Roma, Italy.\linebreak
molle@mat.uniroma2.it
}

\end{center}

\vspace{5mm}



{\small {\sc \noindent \ \ Abstract} -
The paper concerns with positive solutions of problems of the type
$-\Delta u+a(x)\, u=u^{p-1}+\varepsilon u^{2^*-1}$ in
$\Omega\subseteq\R^N$, $N\ge 3$, $2^*={2N\over N-2}$, $2<p<2^*$.
Here $\Omega$ can be an exterior domain, i.e. $\R^N\setminus\Omega$
 is bounded, or the whole of $\R^N$.
The potential $a\in L^{N/2}_{\loc}(\R^N)$ is assumed to be strictly positive and
such that there exists $\lim_{|x|\to\infty}a(x):=a_\infty>0$.
First, some existence results of ground state solutions are proved.
Then the case $a(x)\ge a_\infty$ is considered, with $a(x)\not\equiv
a_\infty$ or $\Omega\neq\R^N$.
In such a case, no ground state solution exists and the existence
of a bound state solution is proved, for small $\e$.

\vspace{3mm}



 \vspace{1mm}


}


\sezione{Introduction and main results}


This paper deals with a class of problems of the type
$$
(P_{\varepsilon})\quad
\left\{
\begin{array}{ll}
-\Delta u+a(x)\, u=u^{p-1}+\varepsilon u^{2^*-1} & \mbox{in $\Omega$}, \\
\noalign{\medskip}
u>0 & \mbox{in $\Omega$}, \\
\noalign{\medskip}
u\in H^1_0(\Omega) \\
\end{array}
\right.
$$
where $\Omega\subseteq\R^N$, $N\ge 3$,  and we consider both the case
$\Omega=\R^N$ and $\R^N\setminus \Omega$ bounded with smooth boundary;
$2^*=\frac{2N}{N-2}$ is the critical Sobolev exponent,
$\varepsilon>0$, $2<p<2^*$ and on the potential we assume
\begin{equation}
\label{Ha}
a\in L^{N/2}_{\loc}(\R^N),\qquad
\lim_{|x|\to\infty}a(x)=a_\infty,\qquad a(x)\ge a_0>0\ \mbox{ a.e. in }\R^N.
\end{equation}

Problem $(P_\e)$ has a variational structure: its solutions
correspond to the nonnegative functions that are critical points of
the functional $E_\e:H^1_0(\Omega)\to\R$ defined by
$$
E_{\e}(u)=\frac{1}{2}\int_{\Omega} (|\nabla u|^2+a(x)\, u^2)\,dx
-\frac{1}{p}\int_{\Omega} |u|^p\,dx -\frac{\e}{2^*}\int_{\Omega}
|u|^{2^*}\,dx.
$$

Problems of the type $(P_\e)$ have been widely studied: it is well
known that they come from problems in Physics and  Mathematical
Physics like Schr\"odinger equations and Klein-Gordon equations,
and from other applied and theoretical sciences.
From a mathematical point of view, problems like $(P_\e)$ present a
number of difficulties related to the lack of compactness due both to
the critical exponent and to the unboundedness of the domain.
If $\R^N\setminus \Omega$ is a ball and $a$ is radially symmetric,
then a classical feature is to employ the compactness of the embedding of
$H^1_{\rad}(\R^N)\hookrightarrow L^p(\R^N)$, that allows to recover
existence results and qualitative properties of solutions for
equations of the type $-\Delta u+a(|x|)u=f(|x|,u)$
(\cite{BereLionsI-II,S}).

For exterior domains and potentials without any symmetry, several
papers treat the subcritical case, i.e. $\e=0$ in $(P_\e)$, starting
from the seminal papers \cite{BC}, concerning the autonomous case,
and \cite{BLi,BPLL}, concerning also the nonautonomous case;
in those papers the authors analyze how the lack of compactness works.
Then, many papers deal with the non autonomous case, in the subcritical
setting (see \cite{CM03,CMP07,MP98,MoM} and references therein).

When $\e>0$ it is interesting to study problem $(P_\e)$ because there
is an overlapping between the effects of the subcritical and the
critical growth in the nonlinearity.
Actually, if $\e>0$, the analysis of the Palais-Smale sequences done in the
subcritical case does not work, so that it is not possible to apply in
a straight way the methods developed in the cited papers.
Indeed, some concentration phenomena can appear, related to the
critical nonlinearity.
Of course, if $\e$ is very large the effect of the critical nonlinearity is
relevant, as one can see, for example, in \cite{MPis03}.
In  \cite{MPis03} the authors prove the existence of solutions of
problems similar to $(P_\e)$, in bounded
domains, and point out some concentration effects as $\e\to\infty$.

Here we want to analyze the problem for small $\e$, so that we have a
critical perturbation of the subcritical case.
Then, besides the analysis of the lack of compactness as in \cite{BC},
we make a further study of the Palais-Smale sequences, that takes into
account the concentration phenomena in the spirit of
\cite{BC90,CFS,PLL85,Struwe84}.
We emphasize that in the problems considered in this paper, in order
to study the compactness, it is not possible to use only the classical
analysis of the compactness developed in the subcritical case, nor the
classical analysis developed in the critical case, but some delicate
estimates involving both cases need (see Proposition \ref{PSmin}).
The aim of this analysis will be not only to show that compactness is
restored below a ``bad energy level'', but also that it is restored in a
suitable range above this ``bad level''.
This done, we can recover a result similar to the well known result
stated in \cite{BPLL} about the existence of a bound state solution in
the subcritical case.

\vspace{2mm}

The first results we prove concern ground state solutions.

\begin{teo}
\label{Tgs1}
If $\Omega=\R^N$, $a(x)$ verifies (\ref{Ha}) and
\beq
\label{Hgs}
a(x)\le a_\infty\ \mbox{ a.e. in }\R^N,
\eeq
then there exists $\e_0>0$ such that problem $(P_\e)$ has a ground state solution for every $\e\in(0,\e_0)$.
\end{teo}

In Theorem \ref{Tgs1}  assumption (\ref{Hgs}) allows to apply in a straight way
concentration-compactness arguments.
Let us consider now the case in which at least one of the assumptions
of Theorem \ref{Tgs1} is not true, that is either $a(x)>a_\infty$ in a
positive measure subset of $\R^N$, or $\Omega\neq\R^N$.
Then, the existence of a ground state solution is not guarantee.
To check the existence of such a solution, the potential
$a(x)$ has to be below $a_\infty$ in a suitable large
region of $\Omega$, to balance the effect of the boundary of $\Omega$
or of the part of $\R^N$ in which $a(x)$ has higher values than $a_\infty$.
In order to state a quantitative assumption, we introduce the problem
$$ (P_{\infty})\quad
\left\{
\begin{array}{ll}
-\Delta u+a_{\infty}u=|u|^{p-2}u& \mbox{in $\R^N$}, \\
\noalign{\medskip}
u\in H^1(\R^N).
\end{array}
\right.
$$
Then, we denote by $w$ the ground state, positive, radial solution of
$(P_\infty)$ and we call
\beq
\label{w_z}
w_z(x):=\vartheta(x)\, w(x-z),\qquad z\in\R^N,
\eeq
where $\vartheta\equiv 1$ if $\Omega=\R^N$, otherwise $\vartheta$  is a cut-off function
verifying
\beq
\label{etheta}
\vartheta\in
\cC^\infty(\R^N,[0,1]),\qquad
\left\{\begin{array}{ll}
\vartheta(x)=1&\mbox{ if }\dist(x,\R^N\setminus\Omega)\ge 1\\
\vartheta(x)=0&\mbox{ if }x\in\R^N\setminus\Omega.
\end{array}\right.
\eeq

\begin{teo}
\label{Tgs2}
Assume that  $a(x)$ verifies (\ref{Ha}).
If there exists $z\in\R^N$ such that
\beq
\label{18.24}
\int_\Omega\left|\D\left({w_z\over \|w_z\|_{L^p(\Omega)}}\right)\right|^2\,
  dx+\int_\Omega a(x)\left({w_z\over\|w_z\|_{L^p(\Omega)}}\right)^2\,
  dx<
{\|w\|_{H^1(\R^N)}^2\over \|w\|_{L^p(\R^N)}^2},
\eeq
then problem $(P_\e)$ has a ground state solution for small $\e$.
\end{teo}

We point out that the r.h.s. in (\ref{18.24}) is a constant
independent of the domain and the potential $a(x)$, and we observe that if
$a(x)\equiv a_\infty$ and $\Omega=\R^N$ then in
(\ref{18.24}) the equality holds for every $z$ in $\R^N$.

Consider now $a(x)\ge a_\infty$.
In Proposition \ref{P4.1} we state that if $\Omega\neq\R^N$ or $a(x)\neq
a_\infty$,  then a ground state solution for $(P_\e)$ does not exist.
In this setting, to find a solution one has to look at higher energy
critical levels and this is more difficult than the minimizing
problem.
A first difficulty to be faced concerns compactness above the ground
state of some related limit problems.
By the concentration phenomenon due to the critical nonlinearity, the
problem to be considered is
$$
(CP_\e)\quad\left\{
\begin{array}{ll}
-\Delta u=\e |u|^{2^*-2}u & \mbox{in $\R^N$}, \\
\noalign{\medskip}
\lim_{|x|\to\infty}u(x)=0,
\end{array}
\right.
$$
while the natural limit problem related to the translations is
$$
(P_{\e,\infty})\quad
\left\{
\begin{array}{ll}
-\Delta u+a_{\infty}u=|u|^{p-2}u+\e |u|^{2^*-2}u & \mbox{in $\R^N$}, \\
\noalign{\medskip}
u\in H^1(\R^N). \\
\end{array}
\right.
$$
A proof of the existence of a ground state solution of $
(P_{\e,\infty})$, for small $\e$, is proved in \cite{AF}.
Looking for least energy solutions of $(P_{\e,\infty})$, the minimization problem
to deal with is
\beq
\label{912}
m_\e:=\inf_{\cN_{\e,\infty}}E_{\e,\infty},
\eeq
where
$$
E_{\e,\infty}(u)=\frac{1}{2}\int_{\R^N} (|\nabla u|^2+a_{\infty}u^2)\,dx
-\frac{1}{p}\int_{\R^N} |u|^p\,dx -\frac{\e}{2^*}\int_{\R^N}
|u|^{2^*}\,dx
$$
and
$$
\cN_{\e,\infty}=\left\{u\in H^1(\R^N)\setminus\{0\}:\enskip
  E_{\e,\infty}'(u)[u]=0\right\}.
$$
Testing the functional $E_{\e,\infty}$ on a concentrating sequence of
least energy solutions of $(CP_\e)$,  in Proposition \ref{Rcrit} we show that
for all $\e>0$
\beq
\label{1106}
m_\e\le {1\over N}\, S^{N/2}\,
\left({1\over\e}\right)^{N-2\over 2},
\eeq
where $S$ is the best Sobolev constant.
We observe that the value ${1\over N}\,
S^{N/2}\,\left({1\over\e}\right)^{N-2\over 2}$ in (\ref{1106}) is the
ground state level of the solutions of problem $(CP_\e)$, that is
$$
{1\over N}\, S^{N/2}\,
\left({1\over\e}\right)^{N-2\over 2}=\min\bigg\{{1\over 2}\int_{\R^N}|\D u|^2dx-{\e\over
    2^*}\int_{\R^N}|u|^{2^*}dx\ :\ u\in \cD^{1,2}(\R^N),\
$$
\beq
\label{1824}
\phantom{********************}
\int_{\R^N}|\D u|^2dx=\e \int_{\R^N}|u|^{2^*}dx\bigg\}
\eeq
(see Proposition \ref{Rcrit}).
So, the compactness cannot hold at the level ${1\over N}\, S^{N/2}\,
\left({1\over\e}\right)^{N-2\over 2}$ neither for problem $(P_\e)$
nor for problem $(P_{\e,\infty})$, by the concentration phenomenon described.
Here  we give an alternative proof of the existence of solutions of
(\ref{912}), that will be useful in the
paper.
As a consequence of that proof, we get that actually $m_\e<{1\over
  N}\, S^{N/2}\, \left({1\over\e}\right)^{N-2\over 2}$, for small $\e$ (see Theorem
\ref{Tm_e} and Corollary \ref{R2.4}).
Nevertheless, neither unicity nor nondegeneracy of the positive
solution of $(P_{\e,\infty})$ are known.
Hence, it is not possible to obtain a complete picture of the lack of
compactness, as in the purely subcritical or critical case.
Anyway, a local Palais-Smale condition can be restored for small $\e$
by using the solutions of $(P_\infty)$.
This done, we can prove the existence of a solution both for the
autonomous and for the non autonomous problem, in $\R^N$ or in
exterior domains.

In \cite{AF} the authors consider problem $(P_\e)$ in the autonomous
case $a(x)\equiv a_\infty$ and they found a solution assuming that $\R^N\setminus
\Omega$ is contained in a small ball.
That result here is improved, because we have no assumption
on the size of $\R^N\setminus \Omega$.
In order to find a solution for every exterior domain, a fundamental
tool is a fine estimate of the interactions of ``almost minimizing''
functions.
Indeed, this estimate allows us to work in a suitable compactness
range (see Lemma \ref{lem:A_varepsilon<2m}).

Our result is the following
\begin{teo}
\label{T}
Assume that $a(x)$ verifies (\ref{Ha}) and
\beq
\label{Has}
a(x)\ge a_\infty,
\qquad\int_{\R^N}(a (x)-a_\infty)|x|^{N-1}e^{2\sqrt{a_\infty}|x|}\,dx<\infty,
\eeq
then there exists $\widehat{\e}>0$ such that for any $0<\e<\widehat{\e}$ problem
$(P_{\e})$ has at least one positive solution, that is a bound state
solution when $\Omega\neq\R^N$ or $a(x)\not\equiv a_\infty$.
\end{teo}

\begin{rem}
\label{R}
{\em
If both $\Omega=\R^N$ and  $a\equiv a_\infty$ hold, problem $(P_\e)$
is nothing but $(P_{\e,\infty})$ and  Theorem \ref{T} coincides with
Theorem \ref{Tm_e}.
}\end{rem}

The paper is organized as follows: in Section 2 we introduce some
notations and recall some known facts we use; Section 3 deals with
ground state solutions; in Section 4 the proof of Theorem \ref{T} is
developed, moreover we report some remarks that describe the asymptotic
shape of the solution given by Theorem \ref{T} and a way to use it to
get multiplicity results (see Remarks \ref{Rf1} and \ref{Rf2}).


\sezione{Notations and preliminary results}


Without any loss of generality we may assume $a_\infty=1$, up to a
rescaling, and $0\in\R^N\setminus\Omega$ if $\Omega\ne\R^N$.
Throughout the paper we make use of the following notation:

{\small
 \begin{itemize}
\item $H^{1}(\R^{N})$ is the usual Sobolev space endowed with
  the standard scalar product and norm
\begin{displaymath}
(u, v):=\int_{\R^N}(\nabla u \nabla v+ uv)dx;\qquad
\|u\|^{2}:=\int_{\R^N}\left(|\nabla u|^{2}+ u^{2}\right)dx.
\end{displaymath}
We shall use also the equivalent norm
$$
\|u\|_a^2:=\int_{\Omega}\left(|\nabla u|^{2}+ a(x)u^{2}\right)dx.
$$
\item $H^{-1}$ denotes the dual space of $H^1(\R^N)$.
\item $\cD^{1,2}(\R^N)$ is the closure of $\cC^{\infty}_0(\R^N)$ with
  respect to the norm $\|u\|_{\cD}:=\left(\int_{\R^N}|\D
    u|^2dx\right)^{1\over 2}$.
\item $L^q(\mathcal{O})$, $1\leq q \leq \infty$, $\mathcal{O}\subseteq
  \R^N$ a measurable set,
  denotes the Lebesgue space, the norm in $L^q(\mathcal{O})$ is denoted
  by $|\cdot|_{L^q( \mathcal{O})}$ when $\mathcal{O}$ is a proper measurable
  subset of $\R^N$ and by $|\cdot|_q$ when $\mathcal{O}=\R^N$.
\item For $u\in H^1_0(\Omega)$ we denote by $u$ also the function in
  $H^1(\R^N)$ obtained setting $u\equiv 0$ in $\R^N\setminus\Omega$.
\item $S$ denotes the best Sobolev constant, namely
$$
S=\min_{u\in\cD^{1,2}(\R^N)\setminus\{0\}}{\int_{\R^N}|\D  u|^2dx
\over
|u|^{2}_{2^*}}.
$$
\item For any $\rho>0$ and for any $z\in \R^N$, $B_\rho(z)$
  denotes the ball of radius $\rho$ centered at $z,$ and for any
  measurable set $ \mathcal{O} \subset \R^N, \ |\mathcal{O}|$ denotes
  its Lebesgue measure.
\item $c,c', C, C', C_i,\ldots$ denote various positive constants.

\end{itemize}
}

When $\e=0$, $(P_\e)$ becomes
$$
(P)\quad
\left\{
\begin{array}{ll}
-\Delta u+a(x)\, u=u^{p-1} & \mbox{in $\Omega$}, \\
\noalign{\medskip}
u>0 & \mbox{in $\Omega$}, \\
\noalign{\medskip}
u\in H^1_0(\Omega)\\
\end{array}
\right.
$$
and the related action functional is $E:H^1_0(\Omega)\to\R$ defined by
$$
E(u)=\frac{1}{2}\int_{\Omega} (|\nabla u|^2+a(x)\, u^2)\,dx
-\frac{1}{p}\int_{\Omega} |u|^p\,dx.
$$
Furthermore, we denote by $E_{\infty},E_{\e,\infty}:H^1(\R^N)\to\R$
the functionals related to $(P_{\infty})$ and
$(P_{\e,\infty})$ respectively, defined by
$$
E_{\infty}(u)=\frac{1}{2}\int_{\R^N} (|\nabla u|^2+u^2)\,dx
-\frac{1}{p}\int_{\R^N} |u|^p\,dx,
$$
$$
E_{\e,\infty}(u)=\frac{1}{2}\int_{\R^N} (|\nabla u|^2+u^2)\,dx
-\frac{1}{p}\int_{\R^N} |u|^p\,dx -\frac{\e}{2^*}\int_{\R^N}
|u|^{2^*}\,dx.
$$
In a standard way, we consider the following Nehari manifolds:
$$
\cN=\left\{u\in H^1_0(\Omega)\setminus\{0\}:\enskip E'(u)[u]=0\right\},
$$
$$
\cN_{\e}=\left\{u\in H^1_0(\Omega)\setminus\{0\}:\enskip E_{\e}'(u)[u]=0\right\},
$$
$$
\cN_{\infty}=\left\{u\in H^1(\R^N)\setminus\{0\}:\enskip E'_\infty(u)[u]=0\right\},
$$
$$
\cN_{\e,\infty}=\left\{u\in H^1(\R^N)\setminus\{0\}:\enskip
  E_{\e,\infty}'(u)[u]=0\right\}.
$$
Remark that there exists $c>0$ independent of small
$\e$ such that
\beq
\label{eN}
\|u\|\ge c\quad\forall u\in \cN_{\e,\infty},
\qquad
\|u\|_a\ge c\quad\forall u\in \cN_{\e},
\eeq
indeed
$$
0=\|u\|^2-|u|_p^p-\e|u|_{2^*}^{2^*}\ge\|u\|^2-c_1\|u\|^p-c_1\e\|u\|^{2^*},
\quad\forall u\in \cN_{\e,\infty},
$$
$$
0=\|u\|_a^2-|u|_p^p-\e|u|_{2^*}^{2^*}\ge\|u\|_a^2-c_2\|u\|^p_a-c_2\e\|u\|_a^{2^*},
\quad\forall u\in \cN_{\e}.
$$
Straight computations allow to state the following

\begin{lemma}
\label{Lnehari}
Let $u\in H^1_0(\Omega)\setminus\{0\}$ and $v\in H^1(\R^N)\setminus\{0\}$,
then:
\begin{itemize}
\item $tu\in \cN$ if and only if
$\displaystyle{t=\left(\frac{\|u\|_a^2}{|u|_p^p}\right)^{1\over p-2}}$;
\item
$tv\in \cN_{\infty}$ if and only if
$\displaystyle{t=\left(\frac{\|v\|^2}{|v|_p^p}\right)^{1\over p-2}}$;
\item $tu\in \cN_{\e}$ if and only if
$\|u\|_a^2=t^{p-2}|u|_p^p+\e t^{2^*-2}|u|_{2^*}^{2^*}$;
\item
$tv\in \cN_{\e,\infty}$ if and only if
$\|v\|^2=t^{p-2}|v|_p^p+\e t^{2^*-2}|v|_{2^*}^{2^*}$.
\end{itemize}
Moreover, $t_u>0$ such that $t_u u\in \cN$ is characterized as
the unique real value such that
$$
E(t_u u)=\max_{t>0} E(tu)
$$
and  $u\mapsto t_u$ is a continuous map from $H^1_0(\Omega)\setminus\{0\}$
in $\R^+$.
Analogous results hold if we consider $E_{\infty}$, $E_{\e}$ and
$E_{\e,\infty}$ respectively on $\cN_{\infty}$, $\cN_{\e}$ and
$\cN_{\e,\infty}$.
\end{lemma}

Let us define:
\beq
\label{min}
m=\inf_{\cN_{\infty}} E_{\infty},\qquad
m_{\e}=\inf_{\cN_{\e,\infty}} E_{\e,\infty}.
\eeq
We denote by $w$ the unique positive solution, up to translations, of
the problem $(P_{\infty})$;
it is well known that $w\in \cC^{\infty}(\R^N)$, $w$ is radially
symmetric about the origin, and
\begin{equation} \label{eqn:decrescita_w}
w(|x|)\, e^{|x|}\,|x|^{(N-1)/2} \to c \quad \mbox{as $|x|\to +\infty$}, \\
\end{equation}
\begin{equation} \label{eqn:decrescita_w'}
w'(|x|)\, e^{|x|}\,|x|^{(N-1)/2} \to -c \quad \mbox{as $|x|\to +\infty$}, \\
\end{equation}
with $c>0$;
moreover $w\in \cN_{\infty}$ and $E_{\infty}(w)=m$, namely $w$ is the
ground state solution of $(P_{\infty})$ (see
\cite{BereLionsI-II,GNN,K} and also (2.19), (2.20) in \cite{BPLL} for
a precise estimate of $c$ in (\ref{eqn:decrescita_w}) and
(\ref{eqn:decrescita_w'})).

\vspace{2mm}

For the limit problem  $(P_{\e,\infty})$ the following existence result holds.

\begin{teo}
\label{Tm_e}
There exists $\e_0>0$ such that for any $\e\in(0,\e_0)$ problem
$(P_{\e,\infty})$ has a positive radially symmetric ground state
solution $w_\e$.
\end{teo}

\proof
We first observe that $m_\e\le m$, $\forall \e>0$.
Indeed let
$\tau_\e>0$ be such that $\tau_\e w\in\cN_{\e,\infty}$, then
\beq
\label{18.28}
m_\e\le E_{\e,\infty}(\tau_\e w)\le E_\infty(\tau_\e w)\le E_\infty(w)=m.
\eeq

As shown in \cite{S}, by Schwartz symmetrization, in order to solve
the minimization problem
for $m_\e$  we can restrict our considerations to
$$
H^1_r(\R^N)=\{u\in H^1(\R^N)\ :\ u \mbox{ radially symmetric}\},\quad
\cN_r=\cN_{\e,\infty}\cap H^1_r(\R^N).
$$
Let $\{u_n\}_n$ in $\cN_r$ be a minimizing sequence, that is
\beq
\label{e1}
\|u_n\|^2=|u_n|_p^p+\e|u_n|_{2^*}^{2^*},
\eeq
\beq
\label{e1.5}
E_{\e,\infty}(u_n)=\left({1\over 2}-{1\over p}\right)\|u_n\|^2+\left({1\over
    p}-{1\over 2^*}\right)\e|u_n|_{2^*}^{2^*}=m_\e+o(1).
\eeq
Inequalities (\ref{18.28}) and  (\ref{e1.5}) imply that
\beq
\label{956}
\|u_n\|^2\le\left({1\over 2}-{1\over p}\right)^{-1} m_\e+o(1)
\le \left({1\over 2}-{1\over p}\right)^{-1} m+o(1).
\eeq
Observe that from (\ref{e1}), (\ref{eN}), (\ref{956}) and the Sobolev
embedding Theorem it follows the existence of $\e_0>0$ such that, for
all $n\in\N$,
\beq
\label{e2}
|u_n|_p^p\ge\|u_n\|^2-c\, \e \|u_n\|^{2^*}\ge\const>0
\qquad
\forall \e\in(0,\e_0).
\eeq

Now, since $H^1_r(\R^N)$ embeds compactly in $L^p(\R^N)$ (see
\cite{S}) we deduce the existence of $w_\e\in
H^1_r(\R^N)$ such that,
up to a subsequence,
\beq
\label{e3}
u_n\xrightarrow[]{n\rightarrow \infty}
w_\e \quad
\left\{\begin{array}{l}
\mbox{ strongly in }L^p(\R^N)\\
\mbox{ weakly  in }H^1(\R^N)\mbox{ and in }L^{2^*}(\R^N),
\end{array}\right.
\eeq
moreover by (\ref{e2}) $w_\e\neq 0$.
By Ekeland's variational principle the minimizing sequence $\{u_n\}_n$
in $\cN_r$ can be chosen such that
\beq
\label{e4}
E'_{\e,\infty}(u_n)[v]=\lambda_n G'(u_n)[v]+o(1)\|v\|\qquad\forall
v\in H^1_r(\R^N)
\eeq
where, for all $n\in\N$, $\lambda_n\in\R$ is the Lagrange multiplier
and $G(u)=E_{\e,\infty}'(u)[u]$.
By definition of $\cN_{\e,\infty}$,  $G(u_n)=0$ for all $n\in\N$, so
using (\ref{e4}), we deduce
\beq
\label{e4.5}
0=G(u_n)=E_{\e,\infty}'(u_n)[u_n]=\lambda_n G'(u_n)[u_n]+o(1)\|u_n\|.
\eeq
Taking into account that $\|u_n\|$ is bounded and that
$G'(u_n)[u_n]\le c<0$ on $\cN_r$,  we get $\lambda_n=o(1)$.
Then (\ref{e4}) implies that $\forall v\in H^1_r(\R^N)$
\beq
\label{e5}
\int_{\R^N}\D u_n\cdot\D v\,dx+\int_{\R^N}u_nv\,dx
-\int_{\R^N}|u_n|^{p-2}u_nv\, dx
-\e \int_{\R^N}|u_n|^{2^*-2}u_nv\, dx
=o(1)\|v\|.
\eeq
By (\ref{e3}) and (\ref{e5}) the function $w_\e$ verifies
$$
\int_{\R^N}\D w_\e\cdot\D v\,dx+\int_{\R^N}w_\e v\,dx
-\int_{\R^N}|w_\e|^{p-2}w_\e v\, dx
-\e\int_{\R^N}|w_\e|^{2^*-2}w_\e v\, dx
=0,
$$
$\forall v\in H^1_r(\R^N)$, so that by choosing $v=w_\e$ it follows
that $w_\e\in\cN_{\e,\infty}$.

Using again (\ref{e3}), we get
\begin{eqnarray*}
m_\e\le E_{\e,\infty}(w_\e)&=&\left({1\over 2}-{1\over
    2^*}\right)\|w_\e\|^2-\left({1\over p}-{1\over 2^*}\right)|w_\e|_{p}^{p}
\\
&\le&
\liminf_{n\to\infty}\left[\left({1\over 2}-{1\over
      2^*}\right)\|u_n\|^2-\left({1\over p}-{1\over
      2^*}\right)|u_n|_{p}^{p}\right]=m_\e,
\end{eqnarray*}
that is $w_\e$ is the minimizing function we are looking for.
Thus, $w_\e$ solves
\beq
\label{e6}
-\Delta u+u=|u|^{p-2}u+\e |u|^{2^*-2}u\qquad \mbox{ in }\R^N.
\eeq

In order to verify that $w_\e$ is strictly positive we just observe
that $|w_\e|$
too is a minimizer of $E_{\e,\infty}$ constrained on
$\cN_{\e,\infty}$, so we can assume $w_\e\ge 0$.
Furthermore, since $w_\e$ solves (\ref{e6}), $w_\e>0$ as a consequence
of the maximum principle.

\qed

\begin{prop}
\label{Rcrit}
The following estimate holds:
\beq
\label{1106bis}
m_\e\le {1\over N}\, S^{N/2}\,
\left({1\over\e}\right)^{N-2\over 2}
\qquad\forall \e>0.
\eeq
\end{prop}
{\noindent{\bf Sketch of the proof \hspace{2mm}}} 
Since the computations that prove (\ref{1106bis}) are classical, we
only sketch them.
First we verify (\ref{1824}). 
Observe that for every $u\in\cD^{1,2}(\R^N)\setminus\{0\}$ the function
$tu$ verifies 
$\int_{\R^N}|\D(tu)|^2dx=\e\int_{\R^N}|tu|^{2^*}dx$ if and only if
$t=\left[{1\over \e}\, {\int_{\R^N}|\D u|^2dx\over
    \int_{\R^N}|u|^{2^*}dx}\right]^{1\over {2^*-2}}$, and
\begin{eqnarray}
\nonumber
{1\over 2}\int_{\R^N}|\D(tu)|^2dx-{\e\over 2^*}\int_{\R^N}|tu|^{2^*}dx
&=&
\left({1\over 2}- {1\over 2^*}\right)t^2 \int_{\R^N}|\D u|^2dx
\\ \nonumber
&=& {1\over N}\left({1\over \e}\right)^{2\over
  2^*-2}{\left(\int_{\R^N}|\D u|^2dx\right)^{{2\over 2^*-2}+1}\over
\left(\int_{\R^N}|u|^{2^*}dx\right)^{2\over 2^*-2} }
\\ \nonumber
&=& {1\over N}\left({1\over \e}\right)^{N-2\over 2}
\left({\int_{\R^N}|\D u|^2dx
\over
|u|^{2}_{2^*}}\right)^{N\over 2}
\\ \label{1842}
&\ge & {1\over N}\left({1\over \e}\right)^{N-2\over 2} S^{N/2}.
\end{eqnarray}
So we obtain (\ref{1824}) by  (\ref{1842}), taking into account that
the equality in (\ref{1842}) is attained by choosing $u$ as a minimizing
function for the Sobolev constant (see \cite{Au,CGS,T}).

Now, let $\bar U\in\cD^{1,2}(\R^N)$ be a fixed radial function that
realizes the minimum in (\ref{1824}), for example consider $\bar
U(x)={C\over (1+|x|^2)^{N-2\over 2}}$, where $C$ is a normalizing
constant.
In order to prove (\ref{1106bis}), we consider the concentrating
sequence of functions $v_n(x)=\zeta(|x|)\, n^{N-2\over 2}\bar U(nx)$, $n\in\N$,
where $\zeta\in\cC^\infty_0(\R^+,[0,1])$ is a cut-off function such
that $\zeta(s)=1$, for $s\in[0,1]$.
Then we test the functional $E_{\e,\infty}$ on the sequence of functions
$u_n:=t_nv_n$, $n\in\N$, where $t_n$ is such that
$u_n\in\cN_{\e,\infty}$, that is
\beq
\label{1918}
|\D v_n|_2^2+ |v_n|_2^2
=
t_n^{p-2} |v_n|_p^p +\e t_n^{2^*-2} |v_n|_{2^*}^{2^*}.
\eeq
Well known estimates provided in \cite{BN} ensure that
\begin{eqnarray}
\label{ea2}
& & | v_n-n^{N-2\over 2}\bar U(nx)|_{2^*}\xrightarrow[]{n\rightarrow \infty} 0,
\\
\label{eb2}
& &  | \D v_n-\D(n^{N-2\over 2}\bar
U(nx))|_{2}\xrightarrow[]{n\rightarrow \infty} 0,
\\
\label{ec2}
& & \ v_n \xrightarrow[]{n\rightarrow
  \infty} 0,\ \mbox{ in } L^2(\R^N).
\end{eqnarray}
From (\ref{ec2}) and the boundedness of $\{v_n\}_n$ in $L^{2^*}(\R^N)$
we obtain also
$v_n \to 0$ in  $L^p(\R^N)$, by interpolation.
Hence, $t_n\to 1$ follows from (\ref{1918}), -- ,(\ref{ec2}) and so
$$
E_{\e,\infty}(u_n)-\left(
{1\over 2}\int_{\R^N}|\D(n^{N-2\over 2}\bar U(nx) )|^2dx-{\e\over
  2^*}\int_{\R^N}|n^{N-2\over 2}\bar U(nx)|^{2^*}dx\right)
\xrightarrow[]{n\rightarrow \infty} 0.
$$
Since
$$
{1\over 2}\int_{\R^N}|\D(n^{N-2\over 2}\bar U(nx) )|^2dx-{\e\over
  2^*}\int_{\R^N}|n^{N-2\over 2}\bar U(nx)|^{2^*}dx = {1\over
  N}\left({1\over \e}\right)^{N-2\over 2} S^{N/2},\qquad\forall n\in\N,
$$
then (\ref{1106bis}) is proved.

\qed

\begin{cor}
\label{R2.4}
For $\e$ small the following estimate holds:
\beq
\label{935bis}
m_\e<{1\over N}\, S^{N/2}\,
\left({1\over\e}\right)^{N-2\over 2}.
\eeq
\end{cor}

\noindent Indeed, in the proof of Proposition \ref{Rcrit} we have
exhibited a sequence 
$\{v_n\}_n$ of radial functions in $\cN_{\e,\infty}$ that converges weakly to 0 in
$L^{2^*}(\R^N)$ and such that  $E_{\e,\infty}(v_n)\to  {1\over
  N}\left({1\over \e}\right)^{N-2\over 2} S^{N/2}$.
But in the proof of Theorem \ref{Tm_e} we have proved that for small
$\e$ {\em every} minimizing sequence of radial functions converges weakly to
a nonzero minimizing function of $E_{\e,\infty}$ on $\cN_{\e,\infty}$, up to a
subsequence.
Hence (\ref{935bis}) must hold, for small $\e$. 

\vspace{2mm}

Let us give another estimate of $m_\e$,  more precise for small $\e$, and
analyze its  asymptotic behaviour.

\begin{lemma} \label{lem:limit_mu_varepsilon}
For all $\e>0$ the relation  $m_\e\le m$ holds  and
$$
\lim_{\e\to 0} m_{\e}=m.
$$
\end{lemma}

\proof
Inequality  $m_\e\le m$  has been shown in (\ref{18.28}).

Now, for $\e\in(0,\e_0)$ let $w_\e$ be the minimizing function whose existence is stated
in Theorem \ref{Tm_e} and
$t_\e>0$ be such that $t_\e w_\e\in\cN_\infty$, namely
\beq
\label{18.57}
t_\e=\left({\|w_\e\|^2\over |w_\e|_p^p}\right)^{1\over p-2}.
\eeq
Observe that $\|w_\e\|$ is bounded, uniformly with respect to
$\e\in(0,\e_0)$, because 
$$
E_{\e,\infty}(w_\e)=\left({1\over 2}-{1\over
    p}\right)\|w_\e\|^2+\left({1\over p}-{1\over
    2^*}\right)\e|w_\e|_{2^*}^{2^*}=m_\e\le m.
$$
Moreover, $|w_\e|_p^p\ge c>0$ follows from (\ref{e2}).
 As a consequence, $t_\e$ is bounded by (\ref{18.57}) and
\begin{eqnarray}
\nonumber m & \le & E_\infty(t_\e w_\e)=E_{\e,\infty}(t_\e w_\e)+{\e\over 2^*}\int_{\R^N}
(t_\e w_\e)^{2^*}dx\\
\label{658} & \le & E_{\e,\infty}(w_\e)+{\e\over 2^*}\int_{\R^N}
(t_\e w_\e)^{2^*}dx\\
\nonumber & =& m_\e+o(1).
\end{eqnarray}
Inequality (\ref{658}) completes the proof.

\qed


\sezione{Existence of a ground state solution}


In this section we prove Theorems \ref{Tgs1} and \ref{Tgs2}, which
provide some cases in which a least energy solution $\bar u$ of
$(P_\e)$ exists, that is  $\bar u\in\cN_\e$ verifies
$$
E_\e(\bar u)=\min_{\cN_\e}E_\e.
$$
A basic tool to prove the existence of a ground state is the
analysis of the Palais-Smale sequences at a level $c$
((PS$)_c${\em-sequences for short}) below the minimum of the limit
problem $(P_{\e,\infty})$.
We start with the following lemma.

\begin{lemma}
\label{PS>0}
{\em
Let $c\in\R$ and let $\{u_n\}_n$ be a (PS$)_c$-sequence for $E_\e$,
then $\{u_n\}_n$ is bounded and $c\ge 0$.
}
\end{lemma}

\proof
From
$$
E'_\e (u_n)[u_n]=\|u_n\|_a^2-|u_n|_p^p-\e|u_n|_{2^*}^{2^*}=o(1)\|u_n\|
$$
we infer
$$
E_\e (u_n)=\left({1\over 2}-{1\over
    p}\right)\|u_n\|_a^2+\e\left({1\over p}-{1\over
    2^*}\right)|u_n|_{2^*}^{2^*}+o(1)\|u_n\|=c+o(1),
$$
that implies our claims.

\qed

\begin{prop}
\label{PSmin}
Assume that $a(x)$ verifies (\ref{Ha}).
Let $\e>0$ and $\{u_n\}_n$ be a $(PS)_c$-sequence for
$E_\e$ constrained on $\cN_\e$.
If $c<m_\e$ then $\{u_n\}_n$ is relatively compact.
\end{prop}

\no \proof
First, let us observe that the sequence $\{\|u_n\|_a\}_n$ is bounded
away from 0 by (\ref{eN}) and it is bounded above because
$$
E_{\e}(u_n)=\left({1\over 2}-{1\over p}\right)\|u_n\|_a^2+
\left({1\over p}-{1\over  2^*}\right)\e\,
|u_n|^{2^*}_{2^*}\xrightarrow[]{n\rightarrow \infty} c.
$$

Then, arguing exactly as in (\ref{e4}),--,(\ref{e5}), we get that
$\{u_n\}_n$ is a $(PS)_c$-sequence also for the free functional
$E_\e$, namely  $\forall v\in
H^1_0(\Omega)$
\beq
\label{1701}
\int_{\Omega}\D u_n\cdot \D v\,dx+\int_{\Omega} au_nv\,dx
-\int_{\Omega}|u_n|^{p-2}u_nv\, dx
-\e \int_{\Omega}|u_n|^{2^*-2}u_nv\, dx
=o(1)\|v\|.
\eeq

From now on, we denote by $\{u_n\}_n$ not only the sequence
$\{u_n\}_n$ but also its subsequences.

Since $\{u_n\}_n$ is bounded in $H^1_0(\Omega)$, there exists a
function $\bar u\in H^1_0(\Omega)$ such that
\beq
\label{14.46}
u_n\xrightarrow[]{n\rightarrow \infty}\bar
u\qquad\left\{\begin{array}{l}
\mbox{ weakly in }H^1_0(\Omega)\ \mbox{ and in }L^{2^*}(\Omega)\\
\mbox{ strongly in }L^p_{\loc}(\R^N)\ \mbox{ and in }L^2_{\loc}(\R^N)\\
\mbox{ a.e. in }\R^N.
\end{array}\right.
\eeq
By (\ref{14.46}) and (\ref{1701}), $\bar u$ is a weak solution of
$(P_\e)$, hence
\beq
\label{14.47}
\|\bar u\|_a^2=|\bar u|_p^p+\e|\bar u|_{2^*}^{2^*}.
\eeq
We have to prove that $u_n\to\bar u$ in $H^1(\Omega)$.
Assume by contradiction that $u_n\not \to\bar u$ in $H^1(\Omega)$, so
the sequence $v_n:=u_n-\bar u$ verifies $\|v_n\|\ge \hat c>0$,
$\forall n\in\N$.
By (\ref{14.46}) and the Brezis-Lieb Lemma,
\beq
\label{1001}
E_\e(u_n)=E_\e(\bar u)+E_\e(v_n)+o(1)
\eeq
and, since $\bar u$ is a solution of $(P_\e)$,
$\{v_n\}_n$ turns out to be a (PS)-sequence for $E_\e$.
We claim that
\beq
\label{14.55}
|v_n|_{2^*}^{2^*}\ge\tilde c>0.
\eeq
If this is not the case, $u_n\to\bar u$ in $L^{2^*}(\Omega)$ and by
interpolation in $L^p(\Omega)$, because $\{u_n\}_n$ is bounded in
$L^2(\Omega)$.
So, from $\|u_n\|^2_a=|u_n|_p^p+\e|u_n|_{2^*}^{2^*}$ and (\ref{14.47})
we get
$$
\lim_{n\to\infty}\|u_n\|_a^2=|\bar u|_p^p+\e|\bar
u|_{2^*}^{2^*}=\|\bar u\|_a^2
$$
which implies $u_n\to \bar u$ in $H^1(\Omega)$, contradicting our assumption.

Let $\{y_i\}_i=\Z^N$ and let us decompose $\R^N$ in the
$N$-dimensional hypercubes $Q_i$ with unitary sides and vertices
in $y_i$.
Since $v_n\in L^{2^*}(\R^N)$, we can define
$$
d_n=\max_{i\in\N}|v_n|_{L^{2^*}(Q_i)}\qquad\forall n\in\N.
$$
By (\ref{14.55}) and the boundedness of $\{u_n\}_n$ in $H^1(\Omega)$
\begin{eqnarray}
\nonumber
0<\tilde c\le |v_n|_{2^*}^{2^*}& =
&\sum_{i=1}^\infty|v_n|^{2^*}_{L^{2^*}(Q_i)}\\
\label{929}
&\le &d_n^{2^*-2}\sum_{i=1}^\infty|v_n|^2_{L^{2^*}(Q_i)}\le
c\, d_n^{2^*-2}\sum_{i=1}^\infty\|v_n\|^2_{H^1(Q_i)}\\
\nonumber
&\le & c'\, d_n^{2^*-2},
\end{eqnarray}
and, then, $d_n\ge\gamma>0$ $\forall n\in\N$, where $\gamma>0$.

Now, let us call $z_n$ the center of a hypercube $Q_{i_n}$ such that
$$
|v_n|_{L^{2^*}(Q_{i_n})}=d_n
$$
and put
$$
w_n(x)=v_n(x+z_n).
$$
Since $\{v_n\}_n$ is a (PS)-sequence, $\{w_n\}_n$ is a (PS)-sequence,
too.

Setting $Q_0=\left[-{1\over 2},{1\over 2}\right]^N$, one of the
following two cases occurs:
\beq
\label{1657}
\begin{array}{rcl}
\vspace{2mm}
(a)& & \int_{Q_0}|w_n(x)|^pdx\ge c>0\\
(b)& &  \int_{Q_0}|w_n(x)|^pdx\xrightarrow[]{n\rightarrow \infty} 0.
\end{array}
\eeq
Assume first that (\ref{1657})$(a)$ holds.
Then $|z_n|\to\infty$ because $v_n\to 0$ in
$L^p_{\loc}(\R^N)$, so, since $\{u_n\}_n$ is a $(PS)$-sequence and $\bar
u$ is a solution of $(P_\e)$, we get
\beq
\begin{array}{r}
\vspace{2mm}
\displaystyle{\int_{\R^N} \D w_n\cdot \D\phi \, dx +\int_{\R^N} w_n \phi \, dx -
\int_{\R^N} |w_n|^{p-2} w_n\phi\,  dx -
\e \int_{\R^N} |w_n|^{2^*-2} w_n\phi\,  dx=}
\\
\label{15.23}
=\displaystyle{ \int_{\R^N} [1-a(\cdot +z_n)]w_n \phi\,
 dx+o(1)\|\phi\|=o(1)\|\phi\|, \quad\forall \phi\in\cC^\infty_0(\R^N)}.
\end{array}
\eeq
The sequence $\{w_n\}_n$ is bounded in $H^1(\R^N)$, so $\bar w\in
H^1(\R^N)$ exists such that
\beq
\label{16.25}
w_n\xrightarrow[]{n\rightarrow \infty}\bar
w\qquad\left\{\begin{array}{l}
\mbox{ weakly in }H^1(\R^N)\ \mbox{ and in }L^{2^*}(\R^N)\\
\mbox{ strongly in }L^p_{\loc}(\R^N)\ \mbox{ and in }L^2_{\loc}(\R^N)\\
\mbox{ a.e. in }\R^N.
\end{array}\right.
\eeq
Now, from  (\ref{1657})$(a)$, (\ref{15.23}), (\ref{16.25}) we deduce
that $\bar w$ is a nonzero solution of $(P_{\e,\infty})$.
Then $\{w_n-\bar w\}_n$ is a $(PS)$-sequence for $E_{\e,\infty}$ and 
$E_{\e,\infty}(w_n-\bar w)\ge o(1)$ can be deduced arguing as in Lemma
\ref{PS>0}.
Hence, applying the Brezis-Lieb Lemma, we get
$$
c=E_\e(u_n)+o(1)=E_\e(\bar u)+E_{\e,\infty}(\bar
w)+E_{\e,\infty}(w_n-\bar w)+o(1)\ge E_{\e,\infty}(\bar
w)+o(1)\ge m_\e+o(1)
$$
contrary to the assumption $c<m_\e$ and proving that (\ref{1657})
$(a)$ can not be true.

\vspace{2mm}

To conclude the argument, we assume that (\ref{1657}) $(b)$ holds and
show that a contradiction arises again.
Remark that in this case we can also assume that
\beq
\label{1105}
\tilde d_n=\max_{i\in\N}|v_n|_{L^p(Q_i)}=\max_{i\in\N}|w_n|_{L^p(Q_i)}
\xrightarrow[]{n\rightarrow \infty} 0.
\eeq
Indeed, if it is not true, we can argue by substituting $Q_{i_n}$ with a cube
 $Q_{\tilde i_n}$ such that $|v_n|_{L^p(Q_{\tilde i_n})}\ge c_1>0$
 and then proceed as in case  (\ref{1657})$(a)$.
So, let us assume (\ref{1105}).
Then, rewriting the inequalities in (\ref{929}) with the
$L^p$-norm in place of the $L^{2^*}$-norm and $\tilde d_n$ in place of
$d_n$, we obtain
\beq
\label{1010}
|v_n|_p=|w_n|_p \xrightarrow[]{n\rightarrow \infty} 0.
\eeq
Notice that (\ref{1010}) implies
\beq
\label{1647}
w_n \xrightarrow[]{n\rightarrow \infty} 0\qquad\mbox{ in }L^2_{\loc}(\R^N).
\eeq
Now, assume that $\{z_n\}_n$ is bounded, so that in our argument we
can consider $z_n=0$, $\forall n\in\N$. 
Let $R>0$ be such that $|a(x)-1|<\eta$ $\forall x\in\R^N\setminus
B_R(0)$, where $\eta$ is a suitable small constant to be fixed later.
Consider the functionals $\hat f$, $\hat f_\infty:\cD^{1,2}(\R^N)\to\R$ defined by
$$
\hat f(u)={1\over 2}\int_{\R^N}|\D u|^2dx+{1\over
  2}\int_{B_R(0)}(a(x)-1)u^2dx-  {\e\over
  2^*}\int_{\R^N}|u|^{2^*}dx,
$$
$$
\hat f_\infty (u)={1\over 2}\int_{\R^N} |\D u|^2dx- {\e\over
  2^*}\int_{\R^N}|u|^{2^*}dx.
$$
Then, (\ref{1701}), (\ref{1010}), and
(\ref{1647}) imply that $\{w_n\}_n$ is a (PS)-sequence also for $\hat f$.
So, Theorem 2.5 of \cite{BC90} applies: there exist a number
$k\in\N\setminus\{0\}$, $k$ sequences of points $\{y_n^j\}_n$, $1\le
j\le k$, $k$ sequences of positive numbers $\{\sigma_n^j\}_n$, $1\le
j\le k$, with $\sigma_n^j\to 0$ because of (\ref{1647}), such that
$$
w_n(x)=\sum_{j=1}^k(\sigma^j_n)^{-{N-2\over 2}}U_j\left({x-y^j_n\over
    \sigma^j_n}\right)+\vi_n(x),
$$
with $\vi_n\to 0$ in $\cD^{1,2}(\R^N)$ and $U_j$ nontrivial solutions
of
\beq
\label{1009}
-\Delta U(x)=\e|U(x)|^{2^*-2}U(x)\qquad x\in\R^N;
\eeq
moreover,
\beq
\label{1008}
\hat f(w_n)=\sum_{j=1}^k\hat f_\infty(U_j)+o(1).
\eeq
By the estimate of the ground state level of the solutions of
(\ref{1009}) given in (\ref{1824}), we get 
\beq
\label{1011}
\hat f_\infty(U_j)\ge {1\over N} \, S^{N/2}\,\left({1\over \e}\right)^{{N-2\over 2}}.
\eeq
Finally, by (\ref{1001}),  (\ref{1008}),  (\ref{1010}), (\ref{1011})
and Proposition \ref{Rcrit} we have 
\begin{eqnarray}
\nonumber E_\e(u_n) & = & E_\e(\bar u)+E_\e(v_n)+o(1)
\\ \nonumber  &\ge &
E_\e(\bar u)+\hat f(w_n)-{\eta\over 2}|w_n|^2_2-{1\over p}|w_n|^p_p+o(1)
\\ \nonumber
& \ge &
E_\e(\bar u)+\sum_{j=1}^k\hat f_\infty(U_j)-\hat c\eta +o(1)
\\ \label{3.15} &\ge &{1\over N} \, S^{N/2}\,\left({1\over \e}\right)^{{N-2\over 2}}-\hat c\eta +o(1)
\\ \nonumber & \ge  & m_\e-\hat c\eta +o(1)
\\ \nonumber & > & c
\end{eqnarray}
for $\eta$ small and large $n$.
So a contradiction arises because of the assumption $E_\e(u_n)\to c$.

Finally, let us consider $|z_n|\to\infty$.
In such a case, the argument developed in the case $\{z_n\}_n$ bounded
can be repeated in an easier way. 
Indeed by (\ref{14.46}), (\ref{1010}) and  (\ref{1647}) we can simply 
consider the functional $\hat f_\infty$ in place of $\hat f$ and get
a contradiction  with $E_\e(u_n)\to c<m_\e$ as in (\ref{3.15}). 
So the proof is completed.

\qed

{\noindent{\bf Proof of Theorem \ref{Tgs1}\hspace{2mm}}}
By Remark \ref{R} we may assume that $a(x)\not \equiv 1$.
We claim that
\beq
\label{1119}
\inf_{\cN_\e}E_\e<m_\e.
\eeq
Indeed, let us consider the minimizing function $w_\e$ for
$(P_{\e,\infty})$ introduced in Theorem \ref{Tm_e} and let $t$ be such
that $t w_\e\in\cN_\e$, then
$$
\inf_{\cN_\e}E_\e\le E_\e (t w_\e)<E_{\e,\infty}(t w_\e)\le
E_{\e,\infty}( w_\e)=m_\e.
$$
By (\ref{1119}) and Proposition \ref{PSmin}  the existence of a
minimizing function $\bar u$ for the functional $E_\e$ constrained on
$\cN_\e$ follows.
Arguing as in the proof of Theorem \ref{Tm_e} one can verify that
$\bar u$ is a constant sign function, which can be chosen strictly
positive.

\qed

{\noindent{\bf Proof of Theorem \ref{Tgs2}\hspace{2mm}}}
Let $z\in\R^N$ be such that (\ref{18.24}) holds and $t_\e>0$ be such
that $t_\e w_z\in\cN_\e$.
In order to obtain the statement, it is enough to prove that for small $\e$
\beq
\label{1139}
E_\e(t_\e w_z)<m_\e.
\eeq
Indeed, once (\ref{1139}) is proved,  $\inf_{\cN_\e} E_{\e}<m_\e$ follows and
we can argue as in the proof of Theorem \ref{Tgs1}.

Let $s>0$ be such that $s w_z\in\cN$, namely
$s=\left({\|w_z\|_a^2\over|w_z|_p^p}\right)^{1\over p-2}$.
We claim that (\ref{18.24}) implies
\beq
\label{1659}
E (s w_z)<m.
\eeq
Let us evaluate
\begin{eqnarray}
\nonumber
E(sw_z) & = & \left({1\over 2}-{1\over p}\right)\|sw_z\|_a^2=
\left({1\over 2}-{1\over p}\right)
\left({\|w_z\|_a^2\over|w_z|_p^p}\right)^{2\over p-2}
\|w_z\|_a^2
\\ \label{1700}  & = &
\left({1\over 2}-{1\over p}\right)
\left(\left\|{w_z\over|w_z|_p}\right\|^2_a\right)^{p\over p-2}.
\end{eqnarray}
Observe that, by (\ref{18.24}),
\beq
\label{1701a}
\left\|{w_z\over|w_z|_p}\right\|^2_a<
{\|w\|^2\over|w|_p^2}
\eeq
and that, since $w$ is the ground state of $(P_\infty)$,
\beq
\label{1702}
\|w\|^2=|w|_p^p\quad\mbox{ and }\quad E_\infty(w)=\left({1\over 2}-{1\over p}\right) \|w\|^2=m.
\eeq
Then, putting (\ref{1701a}) in (\ref{1700}) and using (\ref{1702}), we
get (\ref{1659}).

Finally, remark that $t_\e\to s$, as $\e\to 0$, because
$\|w_z\|^2_a=t^{p-2}_\e|w_z|^p_p+\e t^{2^*-2}_\e|w_z|^{2^*-2}_{2^*}$, and that
$m_\e\to m$ as $\e\to 0$, by Lemma \ref{lem:limit_mu_varepsilon}, so
for small $\e$ (\ref{1139}) follows from (\ref{1659}).

\qed


\sezione{Existence of a bound state solution}


In this section we construct the tools for the proof of Theorem
\ref{T} and prove it.
{We assume}  $\Omega\neq \R^N$ or $a(x)\not\equiv 1$.
First we prove that no ground state solution can exist, then we show
that, in spite of the difficulties due to the few information about the
solutions of $(P_{\e,\infty})$, a local compactness can be recovered
in some interval of the functional values.

\begin{prop} \label{P4.1}
Assume $\e\in[0,\e_0)$. Let $a(x)\ge 1$ and suppose that at least one
assumption between $\Omega\neq\R^N$ and $a(x)\not\equiv 1$ holds true, then
\beq
\label{mu=mbis}
\inf_{\cN_{\e}}E_\e=m_{\e}
\eeq
and the minimization problem (\ref{mu=mbis}) has no solution (here we mean
$E_0=E$, $\cN_0=\cN$, $m_0=m$, $\ldots$).
\end{prop}

\proof
Let $u\in \cN_{\e}$ and $t_u\in\R$ be such that $t_u u\in\cN_{\e,\infty}$.
Since $a(x)\geq 1$ a.e. in $\R^N$, we have
$$
m_{\e}\leq E_{\e,\infty}(t_u u)\leq E_{\e}(t_u u)\leq E_{\e}(u).
$$
Hence $\displaystyle{\inf_{\cN_{\e}}E_{\e}\geq m_{\e}}$.
Let us prove that $\displaystyle{\inf_{\cN_{\e}}E_{\e}\leq m_{\e}}$.

First, assume $\Omega\neq\R^N$. 
In order to exhibit a sequence $\{u_n\}_n$ in $\cN_{\e}$ such
that $E_{\e}(u_n)\to m_{\e}$, we define
$u_n=t_n\,[\vartheta(\cdot)\,w_{\e}(\cdot - ne_1)]$, where
$w_{\e}$ is the minimizing function introduced in Theorem \ref{Tm_e},
$e_1$ is the first element of the canonical basis of $\R^N$, $\vartheta$
is the cut-off function introduced in (\ref{etheta})  and
$t_n> 0$ is such that $u_n=t_n\,[\vartheta(\cdot)\,w_{\e}(\cdot -
ne_1)]\in {\cN}_{\e}.$

Let us fix $r>1$ such that $\R^N\setminus\Omega\subset B_{r-1}(0)$,
then
$$
|\vartheta(\cdot)\,w_{\e}(\cdot - ne_1)-w_{\e}(\cdot - ne_1)|_p^p=
\int_{B_{r}(0)} |(\vartheta(x)-1)\,w_{\e}(x - ne_1)|^p\, dx
$$
$$
\leq  \int_{B_{r}(0)} |w_{\e}(x - ne_1)|^p\, dx =
\int_{B_{r}(-ne_1)} |w_{\e}(z)|^p\, dz=
o(1).
$$
Hence
$|\vartheta(\cdot)\,w_{\e}(\cdot - ne_1)|_{p}^p \to |w_{\e}|_p^p$.
Analogously we have
$|\vartheta(\cdot)\,w_{\e}(\cdot - ne_1)|_{2^*}^{2^*} \to |w_{\e}|_{2^*}^{2^*}$.

Likewise, we have
\beq
\begin{array}{rcl}
\|\vartheta(\cdot)\,w_{\e}(\cdot - ne_1)-w_{\e}(\cdot - ne_1)\|_a^2 & \leq
& \hspace{-1mm} \displaystyle{c \int_{B_{r}(0)} \left(|\nabla \,w_{\e}(x - ne_1)|^2+|w_{\e}(x -
  ne_1)|^2\right)\, dx}
\\ \noalign{\medskip}
\label{eqn:theta_1}
& = & \hspace{-1mm}
o(1).
\end{array}
\eeq
Since $a(x)\to 1$, as $|x|\to \infty$, and
$$
\begin{array}{ll}
\|w_{\e}(\cdot - ne_1)\|_a^2 &=\displaystyle{\int_{\R^N} \left[|\nabla w_{\e}(x-ne_1)|^2+a(x)w_{\e}^2(x-ne_1)\right]\,dx=} \\
\noalign{\medskip}
&=\displaystyle{\int_{\R^N} \left[|\nabla w_{\e}(z)|^2+a(z+ne_1)w_{\e}^2(z)\right]\,dz,}\\
\end{array}
$$
we also get $\|w_{\e}(\cdot - ne_1)\|_a^2 \to \|w_{\e}\|^2$ which,
combined with (\ref{eqn:theta_1}), brings to
$\|\vartheta(\cdot)w_{\e}(\cdot - ne_1)\|_a^2 \to \|w_{\e}\|^2$.

\bigskip

Taking into account $u_n=t_n\,[\vartheta(\cdot)\,w_{\e}(\cdot -
ne_1)]\in {\cN}_{\e}$ and Lemma~\ref{Lnehari},
we have
\begin{equation} \label{eqn:t_n}
\|\vartheta(\cdot)w_{\e}(\cdot - ne_1)\|_a^2-
t_n^{p-2}|\vartheta(\cdot)w_{\e}(\cdot - ne_1)|_{p}^p-
\e\,t_n^{2^*-2}|\vartheta(\cdot)w_{\e}(\cdot - ne_1)|_{2^*}^{2^*}=0,
\end{equation}
so that
$$
t_n^{p-2}|\vartheta(\cdot)w_{\e}(\cdot - ne_1)|_{p}^p+
\e\,t_n^{2^*-2}|\vartheta(\cdot)w_{\e}(\cdot - ne_1)|_{2^*}^{2^*}=
\|\vartheta(\cdot)w_{\e}(\cdot - ne_1)\|_a^2=\|w_{\e}\|^2+o(1).
$$
Hence $\{t_n\}_n$ is bounded and, up to a subsequence, $t_n\to t$.
Getting $n\to \infty$ in (\ref{eqn:t_n}) we obtain
$$
\|w_{\e}\|^2-t^{p-2}|w_{\e}|_p^p-\e\,t^{2^*-2}|w_{\e}|_{2^*}^{2^*}=0,
$$
namely $tw_{\e}\in {\cN}_{\e,\infty}$.
Since $w_{\e}\in {\cN}_{\e,\infty}$, we deduce that $t=1$.
It follows that $\|u_n\|^2\to \|w_{\e}\|^2$, $|u_n|_{p}^p \to
|w_{\e}|_p^p$ and
$|u_n|_{2^*}^{2^*} \to |w_{\e}|_{2^*}^{2^*}$.
Then $E_{\e}(u_n)\to E_{\e,\infty}(w_{\e})=m_{\e}$ and we can conclude
$\displaystyle{\inf_{\cN_{\e}}E_{\e}\leq m_{\e}}$.

If $\Omega =\R^N$, then we set $u_n=t_n\, w_{\e}(\cdot - ne_1)$ and the same
argument developed in the case $\Omega\neq \R^N$ shows that
$E(u_n)\to m_\e$, so that $\displaystyle{\inf_{\cN_{\e}}E_{\e}\leq
  m_{\e}}$ holds again.

Now, let us prove that $m_{\e}$ is not attained in $\cN_{\e}$.
By contradiction, assume that $u\in \cN_{\e}$ verifies
$E_{\e}(u)=m_{\e}$.

First, assume $\Omega\neq\R^N$.
Let $t>0$ be such that $t\,u\in\cN_{\e,\infty}$, then
$$
m_\e\le E_{\e,\infty}(tu)\le E_\e(tu)\le E_\e(u)=m_\e,
$$
i.e. $t\, u$ is a minimizing function for $E_{\e,\infty}$ on
$\cN_{\e,\infty}$.
But then the arguments developed in the proof of Theorem \ref{Tm_e}
show that $|t\, u|>0$, contrary to $u\equiv 0$ in
$\R^N\setminus\Omega$.

Now, assume  $\Omega=\R^N$ and $a\not\equiv 1$.
Again, let $t>0$ be such that $tu\in\cN_{\e,\infty}$, then
$$
m_\e=E_\e(u)\ge E_\e(tu)>E_{\e,\infty}(tu)\ge m_\e,
$$
that is a contradiction, and the proof is complete.

\qed

About the compactness, in the subcritical case we remind an almost
classical result (see f.i. \cite{BC}).
\begin{prop}
\label{PS0}
Let $\{v_n\}_n$ be a $(PS)_c$-sequence of $E$, let $c$ belong to the
interval $(m,2m)$, then $\{v_n\}_n$ is relatively compact and, up to a
subsequence, converges to a nonzero function $\bar v\in H^1_0(\Omega)$
such that $E(\bar v)\in(m,2m)$.
\end{prop}

Here we prove:

\begin{prop}
\label{PS}
To  every $\d\in (0,m/2)$ there corresponds $\e_\d>0$ having the
following property:  $\forall \e\in(0,\e_\d)$, $\forall c\in
(m+\d,2m-\d)$, if $\{u_n\}_n$ is a $(PS)_c$-sequence of $E_\e$ constrained on
$\cN_\e$, then $u_n\rightharpoonup \bar u\neq 0$ weakly in
$H^1(\Omega)$.
Moreover $\bar u$ is a critical point of $E_\e$ on $\cN_\e$ and
$E_\e(\bar u)\le c$.
\end{prop}

\no \proof
As in the proof of Proposition \ref{PSmin}, we deduce that
every (PS$)_c$-sequence for the constrained functional is also a
(PS$)_c$-sequence for the free functional, and its weak limit is a
critical point.
Moreover, by Lemma \ref{PS>0} every (PS)-sequence is bounded in
$H^1(\R^N)$, so it has a weak limit in $H^1(\R^N)$.
Arguing by contradiction, then we can assume that there exist
$\bar\d\in(0,m/2)$, a  sequence $\{c_n\}_n$ in $(m+\bar\d,2m-\bar\d)$,
a sequence $\{\e_n\}_n$ in $(0,+\infty)$, with $\e_n\to 0$, and, for every
$n\in\N$, a sequence $\{u^n_k\}_k$ in $H^1_0(\Omega)$ such that
$$
E_{\e_n}(u^n_k)\xrightarrow[]{k\rightarrow \infty} c_n,\qquad
E_{\e_n}'(u^n_k)\xrightarrow[]{k\rightarrow \infty} 0,
$$
$$
u^n_k\xrightarrow[]{k\rightarrow \infty} 0 \quad\mbox{ weakly in }H^1(\Omega).
$$
Since $p$ is subcritical, we can also assume
$$
u^n_k\xrightarrow[]{k\rightarrow \infty} 0 \quad\mbox{  in }L^p_{\loc}(\R^N).
$$
Now, up to a subsequence, $c_n\to\bar c\in
[m+\bar\delta,2m-\bar\delta]$  and by a diagonal argument we build a
sequence $\{v_n\}_n:=\{u^n_{k_n}\}_n$ such that
\beq
\label{12.30}
E_{\e_n}(v_n)\xrightarrow[]{n\rightarrow \infty} \bar c,
\qquad
E_{\e_n}'(v_n)\xrightarrow[]{n\rightarrow \infty} 0,
\qquad v_n\xrightarrow[]{n\rightarrow \infty} 0 \quad\mbox{  in
}L^p_{\loc}(\R^N).
\eeq
Furthermore
\beq
\label{16.17}
E_{\e_n}(v_n)=\left({1\over 2}-{1\over p}\right)\|v_n\|_a^2+
\left({1\over p}-{1\over 2^*}\right){\e_n}\,
|v_n|_{2^*}^{2^*}=
\bar c+o(1)
\eeq
implies $\{\|v_n\|_a\}_n$ bounded.
Hence we obtain
$$
E(v_n)=E_{\e_n}(v_n)+{\e_n\over 2^*}\int_\Omega
|v_n|^{2^*}dx\xrightarrow[]{n\rightarrow \infty} \bar c
$$
$$
\|E'(v_n)\|_{H^{-1}}\le \|E'_{\e_n}(v_n)\|_{H^{-1}}+C\,\e_n\|v_n\|^{2^*-1}
\xrightarrow[]{n\rightarrow \infty} 0,
$$
so that $\{v_n\}_n$ is a $(PS)_{\bar c}$ -sequence of $E$, with $\bar
c\in (m,2m)$.
Then, by Proposition \ref{PS0},  $\bar v\in H^1(\Omega)$,
$\bar v\neq 0$,  exists such that $v_n\to\bar v$, contrary to (\ref{12.30}).

Finally, if $\{u_n\}_n$ is a (PS$)_c$-sequence for $E_\e$, constrained
on $\cN_\e$, and $u_n\rightharpoonup
\bar u$, then $E_\e(\bar u)\le c$ by (\ref{16.17}) with $\e_n\equiv\e$
and $c$ in place of $\bar c$.

\qed


\subsection{Energy estimates}

Here we first construct some test functions to explore some sublevels of the
functional $E_\e$ and we prove some basic estimates on the action of
these test functions.
Later, we introduce a barycenter map to analyse some features of the sublevels.

\vspace{2mm}

Let us set $\Sigma=\partial B_2(e_1)$, where $e_1$ is the first element
of the canonical basis of $\R^N$, and for any $\rho>0$ define the map
$
\psi_\rho: [0,1]\times\Sigma\longrightarrow H^1_0(\Omega)$ by
$$
\psi_{\rho}[s,y](x)= \vartheta(x)\left[(1-s)w(x-\rho e_1)+sw(x-\rho y)\right],
$$
where $w$ is the ground state solution of $(P_{\infty})$ and
$\vartheta$ is the cut-off function defined in (\ref{etheta}).
Let us denote by $t_{\rho,s,y}$ and $\tau_{\rho,s,y}$ the positive
real numbers such that $t_{\rho,s,y}\,\psi_{\rho}[s,y]\in \cN_{\e}$ and
$\tau_{\rho,s,y}\,\psi_{\rho}[s,y]\in \cN$.

\begin{lemma} \label{lem:A_varepsilon<2m}
There exists $\overline{\rho}>0$ and $\cA\in(m,2m)$ such that for any
$\rho>\overline{\rho}$ and for any $\e>0$
$$
{\cal A}_{\e,\rho}=
\max\left\{E_{\e}\left(t_{\rho,s,y}\,\psi_{\rho}[s,y]\right):\enskip s\in [0,1],\,\, y\in\Sigma\right\} < \cA<2m.
$$
\end{lemma}

Before proving Lemma \ref{lem:A_varepsilon<2m}, let us recall two
technical lemmas. 
We refer the readers to \cite{CP2} for the proof of Lemma \ref{L4.4}
while the proof of Lemma \ref{BL} is in  \cite{BLi} (see also Lemma
2.9 in \cite{CM16}). 

\begin{lemma}
\label{L4.4}
For all $a,b\ge 0$, for all $p\ge 2$, the following relation holds true
\[
(a+b)^{p}\ge a^{p}+b^{p}+(p-1)(a^{p-1}b+ab^{p-1}).
\]
\end{lemma}

\begin{lemma}
\label{BL}
If $g \in L^\infty (\R^N)$ and $h\in L^1(\R^N)$ are
such that, for some $\alpha\ge 0$, $b\ge 0$, $\gamma\in \R$
\beq
\label{eBLg}
\lim_{|x|\to \infty} g(x)e^{\alpha |x|}|x|^b=\gamma
\eeq
and
\beq
\label{eBLh}
\int_{\R^N} |h(x)| e ^{\alpha |x|}|x|^bdx<\infty,
\eeq
then, for every $z\in\R^N\setminus\{0\}$,
\[
\lim_{\rho \to \infty} \left(\int_{\R^N}g(x+\rho z)h(x) dx\right)
e^{\alpha |\rho z|} |\rho z|^b
=\gamma \int_{\R^N}h(x) e^{-\alpha\, {x\cdot z\over|z|}}\, dx.
\]
\end{lemma}

{\noindent{\bf Proof of Lemma \ref{lem:A_varepsilon<2m}}}\hspace{2mm}
In this proof we shall consider $r>1$ fixed such that
$\R^N\setminus\Omega\subset B_{r-1}(0)$, if $\Omega\neq\R^N$, and any
fixed $r>1$ if $\Omega=\R^N$.

Let us set $\delta_{\rho}=\left(\rho^{(N-1)/2}\,e^{2\rho}\right)^{-1}$
and, in order to simplify the notations, we omit $s,y$ and
write $t_\rho=t_{\rho,s,y}$, $\tau_\rho=\tau_{\rho,s,y}$ and
$\psi_{\rho}=\psi_{\rho}[s,y]$.
Being  $\tau_{\rho}\,\psi_{\rho}\in \cN$,
$$
\|\tau_{\rho}\,\psi_{\rho}\|_a^2=|\tau_{\rho}\,\psi_{\rho}|_p^p,
\qquad
\tau_{\rho}=\left(\frac{\|\psi_{\rho}\|_a^2}{|\psi_{\rho}|_p^p}\right)^{1/p-2}
$$
hold true, so, for every $\e>0$, we have
\begin{eqnarray}
\nonumber
E_{\e}(t_{\rho}\,\psi_{\rho})
&\leq&  E(t_{\rho}\,\psi_{\rho})\leq
E(\tau_{\rho}\,\psi_{\rho})
\\
\nonumber
& = &
\frac{1}{2}\|\tau_{\rho}\,\psi_{\rho}\|_a^2-\frac{1}{p}|\tau_{\rho}\,\psi_{\rho}|_p^p
\\
\nonumber
&=&\left(\frac{1}{2}-\frac{1}{p}\right)\tau_{\rho}^2\,\|\psi_{\rho}\|_a^2
\\
\label{1134}
& = & \left(\frac{1}{2}-\frac{1}{p}\right)
\left(\frac{\|\psi_{\rho}\|_a^2}{|\psi_{\rho}|_p^2}\right)^{\frac{p}{p-2}}.
\end{eqnarray}
So, to get the statement of the Lemma, we need to estimate the ratio in
the last line of (\ref{1134}).

\vspace{1mm}

\no {\underline {Estimate of $\|\psi_{\rho}\|_a^2$}}:  we have
\beq
\begin{array}{rcl}
\label{1132}
\|\psi_{\rho}\|_a^2&=&\|\vartheta(\cdot)\left[(1-s)w(\cdot-\rho
  e_1)+sw(\cdot-\rho y)\right]\|_a^2
\\
\noalign{\medskip} & \leq & {\displaystyle{\int_{\R^N}}} |\nabla \vartheta(x)|^2\,\Big[(1-s)w(x-\rho e_1)+sw(x-\rho y)\Big]^2\,dx
\\ \noalign{\medskip}
& & +2 {\displaystyle{\int_{\R^N}}} \Big[\vartheta(x) \nabla\vartheta(x)\Big]\cdot
\Big(\left[(1-s)w(x-\rho e_1)+sw(x-\rho y)\right] \cdot
\\  \noalign{\medskip}
& &\hspace{2cm} \cdot \nabla \left[(1-s)w(x-\rho e_1)+sw(x-\rho
  y)\right]\Big)\,dx
\\ \noalign{\medskip}
& & + {\displaystyle{\int_{\R^N}}} \left(\Big|\nabla \left[(1-s)w(x-\rho e_1)+sw(x-\rho
    y)\right]\Big|^2\right.
\\ \noalign{\medskip}
& &\left. \hspace{2cm} +
a(x)\Big[(1-s)w(x-\rho e_1)+sw(x-\rho y)\Big]^2\right)\,dx.
\end{array}
\eeq
Let us evaluate the addends in (\ref{1132}).
By direct computation and since $w$ is a solution of $(P_{\infty})$, we
obtain
$$
\int_{\R^N} \hspace{-2mm}\left(|\nabla\left[(1-s)w(x-\rho
    e_1)+sw(x-\rho y)\right]|^2\hspace{-2mm}+
a(x)[(1-s)w(x-\rho e_1)+sw(x-\rho y)]^2\right)dx
$$
$$
=\left[(1-s)^2+s^2\right]\|w\|^2+
2s(1-s)\int_{\R^N} w^{p-1}(x-\rho e_1) w(x-\rho y)\,dx
$$
\begin{equation} \label{eqn:norma_pp1}
+\int_{\R^N} (a(x)-1)\Big[(1-s)w(x-\rho e_1) +s\, w(x-\rho y)\Big]^2\,dx.
\end{equation}
By Lemma \ref{BL} there exists ${c_1}>0$ such that
\begin{equation} \label{eqn:cerami_molle_2_9}
\lim_{\rho\to\infty} \delta_{\rho}^{-1}\,\int_{\R^N} w^{p-1}(x-\rho
e_1)\, w(x-\rho y)\,dx=
\end{equation}
$$
=\lim_{\rho\to\infty} \delta_{\rho}^{-1}\,\int_{\R^N} w(x-\rho e_1)\,
w^{p-1}(x-\rho y)\,dx={c_1}.
$$
Taking into account assumption (\ref{Has}) and
(\ref{eqn:decrescita_w}), by Lemma \ref{BL} we have
$$
\hspace{-1cm}\int_{\R^N} (a(x)-1)\Big[(1-s)w(x-\rho e_1) +s\, w(x-\rho
y)\Big]^2\,dx
$$
$$
\hspace{1cm}
\leq
{ 2} \int_{\R^N} (a(x)-1)\Big[w^2(x-\rho e_1)+w^2(x-\rho y)\Big]\,dx=o(\delta_{\rho}).
$$
Hence (\ref{eqn:norma_pp1}) becomes
$$
\int_{\R^N} \hspace{-3mm}\left(|\nabla\left[(1-s)w(x-\rho e_1)+sw(x-\rho y)\right]|^2
\hspace{-1mm} +a(x)[(1-s)w(x-\rho e_1)+sw(x-\rho y)]^2\right)\hspace{-1mm}dx
$$
\begin{equation} \label{eqn:norma_pp2}
\leq\left[(1-s)^2+s^2\right]\|w\|^2+2s(1-s)c_1\,\delta_{\rho}+o(\delta_{\rho}).
\end{equation}
Since $\nabla\vartheta$ has support in $B_r(0)$ and $|y|\ge 1$ $\forall
y\in\Sigma$, from (\ref{eqn:decrescita_w}) it follows
\beq
\label{eqn:norma_nablatheta}
\begin{array}{l}
\hspace{-5mm}\displaystyle{\int_{\R^N} |\nabla \vartheta(x)|^2\,\Big[(1-s)w(x-\rho e_1)+sw(x-\rho y)\Big]^2\,dx}
\\
\\
\leq
\displaystyle{2|\D \vartheta|_{\infty}\int_{B_r(0)}\big[w^2(x-\rho
e_1)+w^2(x-\rho y)\big]\, dx =
o(\delta_{\rho})}.
\end{array}
\eeq
Taking into account (\ref{eqn:decrescita_w'}) and arguing as above we obtain
$$
\hspace{-3cm} 2\int_{\R^N} \Big[\vartheta(x) \nabla\vartheta(x)\Big]\cdot
\Big(\left[(1-s)w(x-\rho e_1)+sw(x-\rho y)\right] \cdot
$$
\begin{equation} \label{eqn:norma_theta}
\hspace{3cm} \cdot \nabla \left[(1-s)w(x-\rho e_1)+sw(x-\rho
  y)\right]\Big)\,dx
\end{equation}
$$
={1\over 2} \int_{B_r(0)}
\nabla(\vartheta(x))^2\cdot \nabla \left[(1-s)w(x-\rho e_1)+sw(x-\rho
  y)\right]^2\,dx
= o(\delta_{\rho}).
$$
By (\ref{1132}),  (\ref{eqn:norma_pp2}), (\ref{eqn:norma_nablatheta}) and
(\ref{eqn:norma_theta}) we deduce
\begin{equation} \label{eqn:norma_psi}
\|\psi_{\rho}\|_a^2 \leq \left[(1-s)^2+s^2\right]\|w\|^2+2s(1-s)c_1\,\delta_{\rho}+o(\delta_{\rho}).
\end{equation}

\vspace{1mm}

{\underline {Estimate of $|\psi_{\rho}|_p^p$}}:
since $0\le \vartheta(x)\le 1$ in $\R^N$ and $\vartheta\equiv 1$ in
$\R^N\setminus B_{r}(0)$, we get
$$
\begin{array}{rl}
|\psi_{\rho}|_p^p &=
\displaystyle{\int_{\R^N} \Big|\vartheta(x)\left[(1-s)w(x-\rho
    e_1)+sw(x-\rho y)\right]\Big|^p\,dx} \\ \\
&\displaystyle{\geq \int_{\R^N} \Big|(1-s)w(x-\rho e_1)+sw(x-\rho y)\Big|^p\,dx} \\
\\
&\hspace{5mm} \displaystyle{-\int_{B_{r}(0)} \Big|(1-s)w(x-\rho e_1)+sw(x-\rho y)\Big|^p\,dx.}
\end{array}
$$
By the asymptotic behaviour of $w$,
$$
\hspace{-4cm}\int_{B_{r}(0)} \Big|(1-s)w(x-\rho e_1)+sw(x-\rho y)\Big|^p\,dx
$$
$$
\hspace{4cm}\leq
 2^{p-1} \int_{B_{r}(0)} \left[ w^p(x-\rho e_1)+ w^p(x-\rho y)\right]\,dx
=o(\delta_{\rho}).
$$
Therefore, from Lemma \ref{L4.4} and (\ref{eqn:cerami_molle_2_9}) it follows
\begin{equation}
\begin{array}{rcl}
\label{eqn:norma_psi_p}
\vspace{2mm}
|\psi_{\rho}|_p^p &\geq&
\displaystyle{\int_{\R^N} \Big|(1-s)w(x-\rho e_1)+sw(x-\rho y)\Big|^p\,dx+o(\delta_{\rho})}\\
 &\geq& \left[(1-s)^p+s^p\right]|w|_p^p+(p-1)\left[(1-s)^{p-1}s+(1-s)s^{p-1}\right]c_1\,\delta_{\rho}
+o(\delta_{\rho}).
\end{array}
\end{equation}

\vspace{1mm}

\no {\underline {Estimate of (\ref{1134})}}:
combining estimates (\ref{eqn:norma_psi}) and (\ref{eqn:norma_psi_p})
and taking advantage of a Taylor expansion, we obtain for any $s\in [0,1]$ and $y\in \Sigma$
$$
\frac{\|\psi_{\rho}\|_a^2}{|\psi_{\rho}|_p^2}\leq
\frac{\left[(1-s)^2+s^2\right]\|w\|^2+2s(1-s)c_1\,\delta_{\rho}+o(\delta_{\rho})}
{\Big(\left[(1-s)^p+s^p\right]|w|_p^p+(p-1)\left[(1-s)^{p-1}s+(1-s)s^{p-1}\right]c_1\,\delta_{\rho}
+o(\delta_{\rho})\Big)^{2/p}}
$$
$$
=\frac{(1-s)^2+s^2}{\left[(1-s)^p+s^p\right]^{2/p}}\,\frac{\|w\|^2}{|w|_p^2}+
\gamma(s)\delta_{\rho}+o(\delta_{\rho}),
$$
where
$$
\gamma(s)=\frac{2s(1-s)c_1}{\left[(1-s)^p+s^p\right]^{2/p}\,|w|_p^2}
\left(1-\frac{p-1}{p}\,\frac{(1-s)^2+s^2}{(1-s)^p+s^p}\left[(1-s)^{p-2}+s^{p-2}\right]
\right).
$$
Since $p>2$ we have that $\gamma(1/2)<0$, so there exist $\overline{c}>0$
and a neighbourhood $I(1/2)$ of $1/2$ such that for any $s\in I(1/2)$
and any $y\in\Sigma$
\begin{eqnarray*}
E_{\e}(t_{\rho}\psi_{\rho})&\leq&
\left(\frac{1}{2}-\frac{1}{p}\right)
\left(\frac{\|\psi_{\rho}\|_a^2}{|\psi_{\rho}|_p^2}\right)^{\frac{p}{p-2}}
\\
&\leq & \left(\frac{1}{2}-\frac{1}{p}\right)
\left(\frac{(1-s)^2+s^2}{\left[(1-s)^p+s^p\right]^{2/p}}\,\frac{\|w\|^2}{|w|_p^2}+
\gamma(s)\delta_{\rho}+o(\delta_{\rho})\right)^{\frac{p}{p-2}}
\\
&\leq & \left(\frac{1}{2}-\frac{1}{p}\right)
\left(2^{\frac{p-2}{p}}\,\frac{\|w\|^2}{|w|_p^2}\right)^{\frac{p}{p-2}}
-\overline{c}\delta_{\rho}+o(\delta_{\rho})
\\
&=&2\left(\frac{1}{2}-\frac{1}{p}\right)|w|_p^p-\overline{c}\delta_{\rho}+o(\delta_{\rho})
\\
&=&
2m-\overline{c}\delta_{\rho}+o(\delta_{\rho}),
\end{eqnarray*}
where we have used $\|w\|^2=|w|_p^p$ and
$E_{\infty}(w)=\left(\frac{1}{2}-\frac{1}{p}\right)|w|_p^p=m$.

Similar computations show that for any $s\in [0,1]\setminus I(1/2)$
and $y\in\Sigma$ we have
$$
\lim_{\rho\to \infty}
\max\{E_{\e}(t_{\rho}\psi_{\rho}):\enskip s\in [0,1]\setminus
I(1/2),\,\, y\in\Sigma\}
$$
$$
\leq \max\left\{
\left(\frac{(1-s)^2+s^2}{\left[(1-s)^p+s^p\right]^{2/p}}\right)^{\frac{p}{p-2}}\,
m:\enskip s\in [0,1]\setminus I(1/2)\right\}<2m.
$$

Finally, we may conclude that the relation
$$
{\cal A}_{\e,\rho}=
\max\left\{E_{\e}\left(t_{\rho,s,y}\,\psi_{\rho}[s,y]\right):\enskip
  s\in [0,1],\,\, y\in\Sigma\right\} < 2m
$$
holds true for $\rho$ large enough, independent of $\e>0$.

\qed

\begin{cor} \label{cor:A_varepsilon<2m}
There exist $\overline{\rho},\overline{\e}>0$ such that for any $\rho>\overline{\rho}$
and for any $\e\in (0,\overline{\e})$
$$
{\cal A}_{\e,\rho}=
\max\left\{E_{\e}\left(t_{\rho,s,y}\,\psi_{\rho}[s,y]\right):\enskip s\in [0,1],\,\, y\in\Sigma\right\} < 2m_\e.
$$
\end{cor}

\proof
It is a direct consequence of Lemmas
\ref{lem:limit_mu_varepsilon} and \ref{lem:A_varepsilon<2m}.

\qed

\vspace{3mm}

The following {\em {definition of barycenter}} of a function $u \in
H^1(\R^N)\setminus \left\lbrace 0 \right\rbrace, $ has been
introduced in
\cite{CP03}.
We set
\beq
\label{mu}
\mu(u)(x)=\frac{1}{|B_1(0)|}\int_{B_1(x)}|u(y)|dy\qquad\qquad x\in\R^N
\eeq
and we remark that  $\mu(u)$ is bounded and continuous, so  we can
introduce the function
\beq
\label{hat}
\hat{u}(x)=\left[\mu(u)(x)-\frac{1}{2}\max \mu(u)\right]^{+}\qquad\qquad x\in\R^N,
\eeq
that is continuous and has compact support.
Thus, we can set
$\beta: H^1(\R^N)\setminus\{0\}\rightarrow \mathbb R^N$
as
$$
\beta(u)=\frac{1}{|\hat{u}|_{1}}\int_{\mathbb R^N}\hat{u}(x)\, x\, dx.
$$
The map $\beta$ has  the following properties:
\begin{eqnarray}
& & \beta \mbox{ is continuous in }H^1(\mathbb R^N)\setminus \{0\};
\label{b1}
\\
&& \mbox{if } u\mbox{ is a radial function, then }\beta(u)=0;
\label{b2}
\\
&& \beta(tu)=\beta(u)
\qquad \forall t\in\R\setminus\{0\},\quad \forall u\in H^1(\mathbb R^N)\setminus
\{0\};
\label{b3}
\\
&& \beta(u(x-z))=\beta(u)+z
\qquad \forall z\in\R^N\quad \forall u\in H^1(\mathbb R^N)\setminus
\{0\}.
\label{b4}
\end{eqnarray}

Let us set
$$
C_0=\inf\{E(u):\enskip u\in \cN,\,\, \beta(u)=0\},
\qquad
C_{0,\e}=\inf\{E_{\e}(u):\enskip u\in \cN_{\e},\,\, \beta(u)=0\}.
$$

\begin{lemma} \label{lem:c_0}
The following facts hold:
\begin{itemize}
\item[$a)$] $C_0>m$;

\item[$b)$] $\displaystyle{\lim_{\e\to 0} C_{0,\e}=C_0}$.
\end{itemize}
\end{lemma}

\proof
Let us prove inequality $a)$.
By Proposition \ref{P4.1}, $C_0\geq m$.
Assume by contradiction that $C_0=m$.
Let $\{{u}_n\}_n$ be a sequence in $\cN$ with $\beta({u}_n)=0$ such that
$E({u}_n)\to m$ and $t_n>0$ be such that $t_nu_n\in
\cN_\infty$, $\forall n\in\N$.
Since $a(x)\geq 1$ a.e. in $\R^N$ we have
\beq
\label{1005}
m\leq E_{\infty}(t_n u_n)\leq E(t_n u_n)\leq E(u_n)=m+o(1),
\eeq
that implies that $\{t_nu_n\}_n$ is a minimizing sequence for
$E_\infty$ on $\cN_\infty$.
Hence there exists a sequence $\{y_n\}_n$ in $\R^N$ such that
$$                                            
t_nu_n(x)=w(x-y_n)+\phi_n(x),\qquad \phi_n\to 0\ \mbox{ strongly in
}H^1(\R^N)
$$
(see \cite[Lemma 3.1]{BC}).
By (\ref{b3}), (\ref{b4}) we have
$$
0=\beta(u_n)=\beta(t_nu_n)=\beta(w(\cdot-y_n)+\phi_n)=
\beta(w+\phi(\cdot+y_n))+y_n.
$$
From $\phi_n\to 0$ strongly in $H^1(\R^N)$ and  (\ref{b1}),
(\ref{b2}), it follows that $\beta(w+\phi(\cdot+y_n))\to \beta(w)=0$,
because $w$ is radially symmetric.
Hence $y_n\to 0$ and $t_nu_n\to w$ strongly in $H^1(\R^N)$.
We shall prove that this is not possible.
If $\Omega\neq\R^N$ then $t_nu_n\equiv 0$ in $\R^N\setminus\Omega$
would imply $w\equiv 0$ in  $\R^N\setminus\Omega$,
contrary to $w>0$ in $\R^N$.
If $\Omega=\R^N$ and $a(x)\not\equiv 1$, then, taking into
account (\ref{1005}), we have
$$
m=E_\infty(w)<E(w)=\lim_{n\to\infty}E(t_nu_n)\le\lim_{n\to\infty}E(u_n)=m,
$$
a contradiction.
So $a)$ is proved.

Let us prove  $b)$.
Let $\e>0$ be fixed and for every $\eta>0$ let  $u_{\eta}\in\cN$ be
such that $\beta(u_{\eta})=0$ and $E(u_{\eta})\leq C_0+\eta$, moreover
let $s_{\eta}>0$ be such that $s_{\eta}u_{\eta}\in\cN_{\e}$.
Then
$$
C_{0,\e}\leq E_{\e}(s_{\eta}u_{\eta})\leq E(s_{\eta}u_{\eta})\leq
E(u_{\eta})\leq C_0+\eta ,
$$
so, by the arbitrary choice of $\eta$, we get
\begin{equation} \label{eqn:co_1}
C_{0,\e}\le C_0\qquad\forall \e>0.
\end{equation}
Let $v_{\e}\in\cN_{\e}$ so that $\beta(v_{\e})=0$ and
$E_{\e}(v_{\e})\leq C_{0,\e}+\e$,
and let $t_{\e}>0$ such that $t_{\e}u_{\e}\in\cN$.
Then
\begin{equation} \label{eqn:co_2}
\begin{array}{rl}
C_0 &\leq \displaystyle{E(t_{\e}v_{\e})=E_{\e}(t_{\e}v_{\e})+\frac{\e}{2^*}|t_{\e}v_{\e}|_{2^*}^{2^*}} \\
\noalign{\medskip}
&\leq \displaystyle{E_{\e}(v_{\e})+\frac{\e}{2^*}|t_{\e}v_{\e}|_{2^*}^{2^*}} \\
\noalign{\medskip}
&\leq \displaystyle{C_{0,\e}+\e+\frac{\e}{2^*}t_{\e}^{2^*}|v_{\e}|_{2^*}^{2^*}.} \\
\noalign{\medskip}
\end{array}
\end{equation}
Now, observe that by (\ref{eqn:co_1})
$$
E_{\e}(v_{\e})=\left(\frac{1}{2}-\frac{1}{p}\right)\|v_{\e}\|_a^2
+\e \left(\frac{1}{p}-\frac{1}{2^*}\right)|v_{\e}|_{2^*}^{2^*}
\leq
C_{0,\e}+\e\leq C_0+\e,
$$
that implies that $|v_{\e}|_{2^*}^{2^*}$ is bounded.
Moreover, taking into account $a(x)\geq 1$ a.e. in $\R^N$ and
$v_{\e}\in\cN_{\e}$, and arguing as
in the proof of (\ref{eN}), we deduce that $\|v_{\e}\|_a\not\to 0$.
Hence, as in (\ref{e2}), we conclude that $|v_{\e}|_p\not\to 0$ too.
Finally, since $t_{\e}v_{\e}\in\cN$, by Lemma \ref{Lnehari} we have that
$\{t_{\e}\}$ is bounded.
So, from (\ref{eqn:co_2}) we infer $\liminf_{\e\to 0}C_{0,\e}\ge C_0$
that, combined with (\ref{eqn:co_1}), gives  $b)$.

\qed

\begin{lemma} \label{lem:c_0_varepsilon}
There exists $\widetilde{\e}>0$ such that for any $\e\in
(0,\widetilde{\e})$ the inequality $C_{0,\e}>{C_0+m\over 2}$ holds.
\end{lemma}

\proof
The assertion follows combining $a)$ and $b)$ of Lemma \ref{lem:c_0}.

\qed

\begin{lemma} \label{lem:c_0_A_varepsilon}
Let ${\cal A}_{\e,\rho}$ be as in Lemma \ref{lem:A_varepsilon<2m}.
Then $\widehat\rho>0$ exists such that $C_{0,\e}\leq {\cal A}_{\e,\rho}$
$\forall \rho>\widehat\rho$, $\forall \e>0$.
\end{lemma}

\proof
We claim that, for $\rho$ large,
$\beta\big(\vartheta(\cdot)w(\cdot-{\rho}y)\big)\cdot y >0$
$\forall y\in\Sigma$.
Indeed, by (\ref{b1})--(\ref{b4}) we have
$$
\Big|\beta\big(\vartheta(\cdot)w(\cdot-\rho y)\big)-\rho y\Big|=
\Big|\beta\big(\vartheta(\cdot+\rho y)\, w\big)\Big|
\xrightarrow[]{\rho\rightarrow \infty}0,
$$
because $\vartheta(\cdot+\rho y)w\to w$ in $H^1(\R^N)$ as $\rho\to \infty$.
Hence
$$
\beta(\vartheta(\cdot)w(\cdot-\rho y))=\rho y+o(1),
$$
that implies the claim.
So, for $\rho$ large, the deformation $\cG:[0,1]\times\Sigma\to
\R^N\setminus\{0\}$ given by
\beq
\label{1552}
\cG(s,y )=s\beta(\psi_{{\rho}}[1,y ])+(1-s)\, y
\eeq
is well defined.
Then, the existence of $({s}_\rho,{y}_\rho)\in [0,1]\times\Sigma$ such that
$\beta(\psi_{{\rho}}[{s}_\rho,{y}_\rho])=0$ follows, because by the continuity
of the maps $\beta$ and $\psi_\rho$ and the invariance of the
topological degree by homotopy we have shown that $0\neq
d(\cG,\Sigma\times[0,1),0)=d(\beta\circ\psi_\rho,\Sigma\times[0,1),0)$.

By (\ref{b3}) we also have
$\beta(t_{{\rho},{s}_\rho,{y}_\rho}\,\psi_{{\rho}}[{s}_\rho,{y}_\rho])$
$=0$.
Since $t_{{\rho},{s}_\rho,{y}_\rho}\,\psi_{{\rho}}[{s}_\rho,{y}_\rho]\in\cN_{\e}$,
the assertion follows.

\qed

\begin{lemma} \label{lem:b_varepsilon}
Let $\tilde \e$ as in Lemma \ref{lem:c_0_varepsilon} and $\e\in (0,\tilde \e)$.
There exists $\widetilde{\rho}>0$ such that for any $\rho>\widetilde{\rho}$
$$
{\cal B}_{\e,\rho}:=
\max\{E_{\e}\left(t_{\rho,1,y}\,\psi_{\rho}[1,y]\right):\enskip y\in \Sigma\}<
C_{0,\e}.
$$
\end{lemma}

\proof
Let us set $t_\rho=t_{\rho,1,y}$ and $\psi_{\rho}=\psi_{\rho}[1,y]$.
By contradiction, assume that there exist $\rho_n\to\infty$ and
$y_n\in\Sigma$ such that
$E_{\e}(t_{\rho_n}\,\psi_{\rho_n})\geq C_{0,\e}$ for every $n\in\N$.

Since $t_{\rho_n}\,\psi_{\rho_n}\in \cN_{\e}$ we can write
\beq
\label{1452}
\begin{array}{rcl}
E_{\e}(t_{\rho_n}\,\psi_{\rho_n})&=&
\displaystyle{\left(\frac{1}{2}-\frac{1}{p}\right)\|t_{\rho_n}\,\psi_{\rho_n}\|_a^2
+\e \left(\frac{1}{p}-\frac{1}{2^*}\right)|t_{\rho_n}\,\psi_{\rho_n}|_{2^*}^{2^*}}
\\
\\
&=&\displaystyle{\left(\frac{1}{2}-\frac{1}{p}\right)t_{\rho_n}^2\|\vartheta
w(\cdot -\rho_n y_n)\|_a^2
+\e \left(\frac{1}{p}-\frac{1}{2^*}\right)t_{\rho_n}^{2^*}|\vartheta
w(\cdot -\rho_n y_n)|_{2^*}^{2^*}}.
\end{array}
\eeq
Observe that in our setting
$0<m\le C_{0,\e}\le E_{\e}(t_{\rho_n}\,\psi_{\rho_n})\leq {\cal
  A}_{\e,\rho}<2m$ and that $0<c\le\|\vartheta
w(\cdot -\rho_n y_n)\|_a\le C<\infty$, $\forall n\in\N$.
Hence from (\ref{1452}) it follows that $0<c_1\le t_{\rho_n} \le C_1<\infty$.
So, up to a subsequence, we can assume $t_{\rho_n} \to t>0$.

Since $\rho_n\to\infty$, the same estimates provided in the proof of
Lemma \ref{lem:A_varepsilon<2m} prove
$E_{\e}(t_{\rho_n}\,\psi_{\rho_n})\to E_{\e,\infty}(tw)$, and we get
\begin{eqnarray*}
C_{0,\e}&\leq &E_{\e,\infty}(tw)=
E_{\infty}(tw)-\frac{\e}{2^*} |tw|_{2^*}^{2^*}
\\
& \leq&
E_{\infty}(w)-\frac{\e}{2^*} |tw|_{2^*}^{2^*}\\
&=&
m-\frac{\e}{2^*} |tw|_{2^*}^{2^*}<m,
\end{eqnarray*}
contrary to Lemma \ref{lem:c_0_varepsilon} and Lemma \ref{lem:c_0} $(a)$.

\qed


\subsection{Proof of Theorem \ref{T}}

Let us recall the values
$$
{\cal A}_{\e,\rho}=
\max\{E_{\e}\left(t_{\rho,s,y}\,\psi_{\rho}[s,y]\right):\enskip s\in
[0,1],\ y\in \Sigma\},
$$
$$
{\cal B}_{\e,\rho}=
\max\{E_{\e}\left(t_{\rho,1,y}\,\psi_{\rho}[1,y]\right):\enskip y\in \Sigma\},
$$
\beq
\label{1644}
C_{0,\e}=\inf\{E_{\e}(u):\enskip u\in \cN_{\e},\,\, \beta(u)=0\}.
\eeq
By Corollary \ref{cor:A_varepsilon<2m} and Lemmas
\ref{lem:A_varepsilon<2m}, \ref{lem:c_0}, \ref{lem:c_0_varepsilon},
\ref{lem:c_0_A_varepsilon} and \ref{lem:b_varepsilon}, the inequalities
\beq
\label{1635}
\left\{
\begin{array}{cl}
(a)&\quad {\cal B}_{\e,\rho}<C_{0,\e}\leq {\cal A}_{\e,\rho}\\
(b)&\quad m<{c_0+m\over 2}<C_{0,\e}\le\cA_{\e,\rho}\le\cA<2m\\
(c)&\quad \cA_{\e,\rho}<2\, m_\e
\end{array}
\right.
\eeq
hold true for every $\rho>\max\{\bar{\rho},\widetilde{\rho},\widehat\rho\}$ and
for every $0<\e<\min\{\bar{\e},\widetilde{\e}\}$.
Let $0<\delta<\min\left\{{m\over 2},2m-{\cal A},{C_0-m\over
    2}\right\}$ and let us consider $\e_\d$ according to Proposition \ref{PS}.

 We  claim that $E_\e$ constrained on $\cN_\e$ has a (PS)-sequence in
$[C_{0,\e},{\cal A}_{\e,\rho}]$ for every
$0<\e<\widehat\e:=\min\{\e_{\delta},\bar{\e},\widetilde{\e}\}$.
This done, the existence of a non-zero critical point $\bar u$ with
$E_\e(\bar u)\le \cA_{\e,\rho}$ follows from Proposition \ref{PS}.

Assume, by contradiction, that no (PS)-sequence exists in
$[C_{0,\e},{\cal A}_{\e,\rho}]$.
Then,  usual deformation arguments imply the existence of $\eta>0$ such that
the sublevel
$E_{\e}^{C_{0,\e}-\eta}:=\{u\in\cN_{\e}:\enskip E_{\e}(u)\leq
C_{0,\e}-\eta\}$ is a deformation retract of the sublevel
$E_{\e}^{{\cal A}_{\e,\rho}}:=\{u\in\cN_{\e}:\enskip E_{\e}(u)\leq {\cal
  A}_{\e,\rho}\}$,
namely there exists a continuous function
$\sigma:E_{\e}^{{\cal A}_{\e,\rho}}\to E_{\e}^{C_{0,\e}-\eta}$ such that
\beq
\label{1553}
\sigma(u)=u\qquad\mbox{ for any }u\in E_{\e}^{C_{0,\e}-\eta}.
\eeq
Furthermore, by (\ref{1635}) $(a)$ we can also assume $\eta$ so small that
\beq
\label{1603}
C_{0,\e}-\eta>{\cal B}_{\e,\rho}.
\eeq
Let us define the map $\cH:[0,1]\times\Sigma\to\R^N$ by
$$
\cH(s,y)=\beta\left(\sigma\big(t_{\rho,s,y}\,\psi_{\rho}[s,y]\big)\right).
$$
By (\ref{1603}),  (\ref{1553}) and by using the map $\cG$ introduced
in (\ref{1552}), we deduce that $\cH$ maps
$\{1\}\times  \Sigma$ in a set homotopically equivalent to
$\rho\Sigma$ (and then to $\Sigma$) in $\R^N\setminus\{0\}$.
Moreover, taking also into account Lemma \ref{Lnehari}, we see
that  $\cH$ is a continuous map.
Hence, by the argument developed in the proof of Lemma
\ref{lem:c_0_A_varepsilon}, a point $(\tilde s,\tilde y)\in [0,1]\times
\Sigma$ must exist, for which 
$$
0=\cH(\tilde s,\tilde y)=\beta(\sigma(t_{\rho,\tilde s,\tilde
  y}\,\psi_{\rho}[\tilde s,\tilde y])).
$$
Then, $E_\e(\sigma(t_{\rho,\tilde s,\tilde
  y}\,\psi_{\rho}[\tilde s,\tilde y]) )\ge C_{0,\e}$, contrary to
$\sigma\big(t_{\rho,s,y}\,\psi_{\rho}[s,y]\big)\in
E_{\e}^{C_{0,\e}-\eta}$ for every $(s,y)\in[0,1]\times \Sigma$, so the
claim must be true.

Let $\bar u\in E_\e^{\cA_{\e,\rho}}$ be the critical point we have found.
To show that $\bar u$ is a constant sign function,
assume, by contradiction, that $\bar u=\bar u^+-\bar u^-$, with $\bar
u^\pm\neq0$.
Multiplying the equation in $(P_\e)$ by $\bar u^\pm$ we deduce that
$\bar u^\pm\in\cN_\e$, so
$$
E_\e(\bar u)=E_\e(\bar u^+)+E_\e(\bar u^-)\ge 2 m_\e,
$$
contrary to (\ref{1635}) $(c)$.

\qed

\begin{rem}
\label{Rf1}
{\em
Let us set
$$
R(\Omega)=\max\{r>0\ :\ \exists x_r\in\R^N\ \mbox{ such that
}B_r(x_r)\subset\R^N\setminus\Omega\}.
$$
Assume $B_{R(\Omega)}(0)\subset \R^N\setminus\Omega$
and call $u_{a,\Omega}$ the solution provided by Theorem
\ref{T}.
Arguing as in \cite{MP98}, the following asymptotic behaviour of
$u_{a,\Omega}$ can be described, as $R(\Omega)\to\infty$, up to some
sequence:
$$
u_{a,\Omega}(x)=w_{1,\e}(x-x_{1,\Omega})+w_{2,\e}(x-x_{2,\Omega})+O(\Omega),
$$
where $O(\Omega)\longrightarrow 0$ in $H^1(\R^N)$, as $R(\Omega)\to\infty$, 
$x_{1,\Omega},x_{2,\Omega}\in\R^N$ verify
$$
|x_{1,\Omega}-x_{2,\Omega}|\longrightarrow \infty\quad
\mbox{ and }\quad {x_{1,\Omega}+x_{2,\Omega}\over 2}\longrightarrow
  0,\qquad\mbox{ as }R(\Omega)\to\infty,
$$
and $w_{1,\e}$, $w_{2.\e}$ are solutions of $(P_{\e,\infty})$.
The same behaviour of $u_{a,\Omega}$ can be obtained considering a
sequence of potentials $a_n(x)$ verifying (\ref{Ha}) and (\ref{Has})
and such that
$$
\lim_{n\to\infty}a_n(x)=\infty\qquad\mbox{ a.e. in }\R^N.
$$
On the contrary, if the capacity of $\R^N\setminus\Omega$ goes to
zero and $|a_n-a_\infty|_{N/2}\to 0$, then $u_{a_n,\Omega}$ converges
to a solution of the limit problem $(P_{\e,\infty})$.
}
\end{rem}

\begin{rem}
\label{Rf2}
{\em
The behaviour of the solution $u_{a,\Omega}$ described in Remark
\ref{Rf1} can be employed to obtain multiplicity of solutions of
$(P_\e)$ when $\Omega=\R^N\setminus\cup_{i=1}^h\omega_i$ and
$a(x)=a_\infty+\sum_{j=1}^k\alpha_j(x)$, with suitable
$\omega_i\subset\hspace{-1mm} \subset\R^N$, $i=1,\ldots,h$, and $\alpha_j\in
L^{N/2}(\R^N)$, $j=1,\ldots,k$.
See  \cite{MoM} for a description of the method.
}\end{rem}

\qed


{\small {\bf Acknowledgements}.
The authors have been supported by the ``Gruppo
Nazionale per l'Analisi Matematica, la Probabilit\`a e le loro
Applicazioni (GNAMPA)'' of the {\em Istituto Nazionale di Alta Matematica
(INdAM)} - Project: Sistemi differenziali ellittici nonlineari
derivanti dallo studio di fenomeni elettromagnetici.

\noindent The first author acknowledges also the MIUR Excellence Department
Project awarded to the Department of Mathematical Sciences, Politecnico of Turin, CUP E11G18000350001.

\noindent The second author acknowledges also the MIUR Excellence Department
Project awarded to the Department of Mathematics, University of Rome
Tor Vergata, CUP E83C18000100006.
}



\end{document}